\documentclass{elsarticle}

\usepackage{lineno,hyperref}
\modulolinenumbers[5]

\journal{%Journal of Differential Equations
}

\usepackage{latexsym}
\usepackage{amsmath}
\usepackage{amssymb}
\usepackage{amsfonts}
\usepackage{verbatim}
\usepackage[latin1]{inputenc}
\usepackage{mathrsfs}
\usepackage{showlabels}
\usepackage{amsfonts}
\usepackage{graphicx}
\usepackage{colortbl,dcolumn}
\usepackage{amsmath}
\usepackage{psfrag}
\usepackage{booktabs}

\newtheorem{tm}{Theorem}[section]
\newtheorem{rk}{Remark}[section]
\newtheorem{ap}{Assumption}[section]
\newtheorem{df}{Definition}[section]
\newtheorem{prop}{Proposition}[section]
\newtheorem{lm}{Lemma}[section]
\newtheorem{cor}{Corollary}[section]

\newtheorem{MP}{Main Problem}[section]

\newcommand{\cc}{\mathbb C}
\newcommand{\ee}{\mathbb E}
\newcommand{\pp}{\mathbb P}
\newcommand{\nn}{\mathbb N}
\newcommand{\rr}{\mathbb R}

\newcommand{\BB}{\mathcal B}
\newcommand{\CC}{\mathcal C}
\newcommand{\LL}{\mathcal L}

\newcommand{\OOO}{\mathscr O}
\newcommand{\FFF}{\mathscr F}
\newcommand{\MMM}{\mathscr M}
\newcommand{\HHH}{\mathscr H}

\begin{document}

\begin{frontmatter}

\title{Well-posedness and Optimal Regularity of Stochastic Evolution Equations with Multiplicative Noises}
\tnotetext[mytitlenote]{This work was funded by National Natural Science Foundation of China (No. 91630312, No. 91530118 and No. 11290142).}

\author[cas]{Jialin Hong}
\ead{hjl@lsec.cc.ac.cn}

\author[HKUST]{Zhihui Liu\corref{cor}}
\ead{liuzhihui@lsec.cc.ac.cn}

\cortext[cor]{Corresponding author.}

\address[cas]{LSEC, ICMSEC, \\
Academy of Mathematics and Systems Science, \\
Chinese Academy of Sciences, \\
 Beijing 100190, China}

\iffalse
\address[HKUST]{Department of Mathematics, \\
The Hong Kong University of Science and Technology, \\
Clear Water Bay, Kowloon, Hong Kong}
\fi

\address[HKUST]{Department of Applied Mathematics, \\
The Hong Kong Polytechnic University, \\
Hung Hom, Kowloon, Hong Kong}

\begin{abstract}
In this paper, we establish the well-posedness and optimal trajectory regularity for the solution of stochastic evolution equations with generalized Lipschitz-type coefficients driven by general multiplicative noises.
To ensure the well-posedness of the problem, the linear operator of the equations is only need to be a generator of a $\CC_0$-semigroup and the proposed noises are quite general, which include space-time white noise and rougher noises.
When the linear operator generates an analytic $\CC_0$-semigroup, we derive the optimal trajectory regularity of the solution through a generalized criterion of factorization method.
\end{abstract}

\begin{keyword}
stochastic evolution equation \sep
multiplicative noise \sep
trajectory regularity \sep
factorization method \sep
Gr\"onwall inequality with singular kernel
\MSC[2010] 
Primary: 60H35 \sep  
Secondary: 60H15
\end{keyword}

\end{frontmatter}

\linenumbers

%\tableofcontents

\section{Introduction}
\label{sec1}

In this paper, we study the well-posedness and optimal trajectory regularity for the solution of the semilinear stochastic evolution equation (SEE)
\begin{align}\label{spde} \tag{SEE}
\begin{split}
{\rm d}X(t)&=(AX(t)+F(X(t))) {\rm d}t+G(X(t)) {\rm d}W(t),
\quad t\in (0,T]; \\
 X(0) &=X_0 
 \end{split}
\end{align}
in a separable Hilbert space $H$, under weak assumptions on the data.
Here $T$ is a fixed positive number and $W:=\{W(t):\ t\in [0,T]\}$ is a ${\bf Q}$-Wiener process with values in another separable Hilbert space $U$ with respect to a stochastic basis $(\Omega,\FFF,(\FFF_t)_{t\in [0,T]},\pp)$.

The well-posedness and the regularity for the solution of an SEE are two fundamental issues in both mathematical and numerical analysis (see, e.g., \cite{CHL17(SINUM)}, \cite{CHJNW17(IMA)}, \cite{CN13(NM)}, \cite{CHL17(JDE)}, \cite{CHLZ17}, \cite{DZ14} and references therein).
These two problems for Eq. \eqref{spde} with finite dimensional multiplicative noises or infinite dimensional affine noises have been studied extensively; see, e.g., 
G. Da prato, S. Kwapie\u{n} $\&$ J. Zabczyk \cite{DKZ87(STO)}, 
N. Krylov \cite{Kry94(PTRF)}, 
S. Tindel, C. Tudor $\&$ F. Viens \cite{TTV04(JFA)} and 
Z. Brzez\'niak, J. van Neerven, M. Veraar $\&$ L. Weis \cite{BNVW08(JDE)} and references therein.
For Eq. \eqref{spde} with finite or infinite dimensional multiplicative smooth noises, we refer to M. Hofmanov{\'a} \cite{Hof13(NoDEA)} and X. Zhang \cite{Zha07(JFA)} where the authors studied conditions on the coefficients and the noises to ensure the existence of a continuous strong solution and the infinitely often differentiability in the spatial variable for the solution of Eq. \eqref{spde}, respectively.
Recently, the authors in \cite{NVW12(SIMA)}, \cite{NVW12(AOP)} and references therein studied the maximal $L^p$-regularity for stochastic convolutions and applied to ${\rm Dom} (-A)^\frac12$-well-posedness of Eq. \eqref{spde} with strong Lipschitz conditions.

One may expect that the solution of Eq. \eqref{spde} with certain assumptions on the initial datum $X_0$, the coefficients $F$ and $G$, inherits the same regularity as the solution of the associated linear SEE
\begin{align}\label{spde-add}
{\rm d}X(t)=AX(t){\rm d}t+{\rm d}W(t), \quad t\in (0,T]; \quad X(0)=X_0.
\end{align}
It is well-known that the unique solution of Eq. \eqref{spde-add} is given by
$X(\cdot)=S(\cdot)X_0+W_A(\cdot)$, where $S(\cdot):=e^{A\cdot}$ is the semigroup generated by $A$ and $W_A(\cdot)$ is the so-called Ornstein--Uhlenbeck process
\begin{align*}
W_A(t):=\int_0^t S(t-r){\rm d}W(r),\quad t\in [0,T].
\end{align*}
If $S(\cdot)$ is a $C_0$-semigroup, then by It\^o isometry $W_A$ defines an $H$-valued stochastic process if and only if
$\int_0^T \|S(r) \|_{\LL_2^0}^2 {\rm d}r<\infty$.

To study the temporal regularity of $W_A$, the authors in \cite{DKZ87(STO)} introduced a factorization formula which was then applied to study numerous SEEs by a lot of researchers in different settings (see, e.g., \cite{HL17(TAMS)}, \cite{NVW08(JFA)} and references therein). 
Under the additional assumption that
$\int_0^T r^{-2\alpha}\|S(r) \|_{\LL_2^0}^2 {\rm d}r<\infty$
holds for some $\alpha\in (0,1/2)$, \cite{DKZ87(STO)} proved that $W_A$ has a continuous version in $H$.
Moreover, if $S$ is supposed to be an analytic $\CC_0$-semigroup satisfying certain properties (see \eqref{ana}), then 
\begin{align} \label{da}
W_A\in \CC^\delta([0,T]; \dot H^\theta)\quad \rm{a.s.}
\end{align} 
for any $\delta,\theta\ge 0$ with $\delta+\theta/2<\alpha$.
The limit case $\alpha=1/2$ is included when ${\bf Q}$ is a trace operator (see \cite[Theorems 5.15 and 5.16]{DZ14}).
Moreover, the optimality of this regularity is shown by a counter-example in \cite{DKZ87(STO)} when $A$ is self-adjoint and positive definite. 
A natural problem whether one can extend the optimal regularity result \eqref{da} to the case $\theta\ge 2\alpha$ is unknown.

Another related interesting problem is to generalize this type of trajectory regularity to the solution of Eq. \eqref{spde} with general data.
An important result is given by A. Jentzen and M. R\"ockner \cite{JR12(JDE)}, where the authors studied the well-posedness and regularity for the solution of Eq. \eqref{spde} driven by a multiplicative trace class noise.
Under the assumptions that $S(\cdot)$ is an analytic $\CC_0$-semigroup, $F:H\rightarrow H$ and 
$G: H\rightarrow \LL_2^0$ are Lipschitz continuous,
$\|G(z)\|_{\LL_2^\gamma}\le C(1+\|z\|_\gamma)$ for some 
$\gamma\in [0,1)$ and any $z\in \dot H^\gamma$, and $X_0\in L^p(\Omega;\dot H^\beta)$ for some 
$\beta\in [\gamma,\gamma+1)$ and $p\ge 2$, they proved the existence of a unique mild solution 
$X\in L^\infty([0,T]; L^p(\Omega; \dot H^\beta))$ such that  
\begin{align}\label{jr}
\ee\Big[\|X(t_1)-X(t_2)\|_\theta^p \Big]
\le C|t_1-t_2|^{\big(\frac12\wedge \frac{\beta-\theta}2\big)p},\quad t_1,t_2\in [0,T],
\end{align}
for any $\theta\in [0,\beta)$ and that $X$ is continuous with respect to $\|\cdot\|_{L^p(\Omega;\dot H^\beta)}$.
It is not clear whether the solution of Eq. \eqref{spde} possesses the trajectory continuity in $\dot H^\beta$.
On the other hand, how to derive the optimal regularity of $X$ for general $\beta$ and $\gamma$ remains open.

As a consequence of \eqref{jr} for $\beta\in [\gamma,\gamma+1)$ and the Kolmogorov continuity theorem,
\begin{align}\label{jr-hol}
\|X(t_1,\omega)-X(t_2,\omega)\|_\theta
\le C(\omega)|t_1-t_2|^\delta,\quad t_1,t_2\in [0,T],
\ \omega\in \Omega,
\end{align}
for any $\delta<[1\wedge (\beta-\theta)]/2-1/p$
and $\theta<\beta-2/p$ provided that $p>2$.
To derive the trajectory continuity of $X$ in $\dot H^\theta$, one needs the restriction that 
$\beta>2/p$ and $\theta<\beta-2/p$.
Indeed, whether $X$ possesses the trajectory continuity in $\dot H^\theta$ when $\beta\le 2/p$ with $\theta\in [0,\beta]$ or $\beta>2/p$ with $\theta\in [\beta-2/p,\beta]$ is still unknown.

The above questions are main motivations for us to study the well-posedness and optimal trajectory regularity for the solution of Eq. \eqref{spde}.
Another motivation is to relax the assumptions on the data $X_0$, $A$, $F$ and $G$ of Eq. \eqref{spde}, which can handle more SEEs in applications.
These motivations lead to the following 

\begin{MP}\label{MP}
To derive the well-posedness and optimal regularity for the solution of Eq. \eqref{spde} under less assumptions on its data.
\end{MP}

To study the well-posedness and optimal trajectory regularity for the solution of Eq. \eqref{spde} and answer the aforementioned questions, we adopt a complete different method compared with \cite{JR12(JDE)}.
It should be noticed that, to establish the well-posedness of Eq. \eqref{spde} under less assumptions on the data, we only need that $S(\cdot)$ is a $C_0$-semigroup. 
To show that the solution is continuous a.s., we need an additional assumption (see Assumption \ref{a3}).
In order to study the trajectory regularity for the solution of Eq. \eqref{spde}, we do not use spectral representation for the linear operator $A$; our main assumption on the operator $A$ is that \eqref{ana} holds.
Thus our well-posedness and continuity results (see Theorems \ref{main1} and \ref{main2}) hold for $C_0$-semigroup and our regularity results hold for analytic $\CC_0$-semigroup (see Theorems \ref{main3} and \ref{main4}).
These results are also new for deterministic evolution equations under our assumptions.
We also mention that the well-posedness and regularity for the solution of Eq. \eqref{spde} in Banach setting have been studied in a companion paper \cite{HHL17}.

The rest of this article is organized as follows.
In the next section, we give our main idea and results and present several concrete examples which satisfy our assumptions.
We prove our well-posedenss as well as trajectory regularity results in Sections \ref{sec3} and \ref{sec4}, respectively.

\section{Main Results}
\label{sec2}

To perform the formulation, let us recall some frequently used notations.
Let $(H, \|\cdot\|_H)$ be a separable Hilbert space and $A: D(A)\subset H\rightarrow H$ be the infinitesimal generator of a $C_0$-semigroup $S(\cdot)$.
In the study of the regularity for the solution of Eq. \eqref{spde}, we assume furthermore that $S(\cdot)$ is an analytic $\CC_0$-semigroup such that the resolvent set of $A$ contains all $\lambda\in \cc$ with $\Re [\lambda]\ge 0$.
Then one can define the fractional powers $(-A)^\gamma$ for $\gamma\in \rr$ of the operator $A$ (see, e.g., \cite[Section 2]{JR12(JDE)} or \cite[Chapter 2.6]{Paz83}).
Let $\dot H^\gamma$ be the domain of $(-A)^\frac\gamma2$ equipped with the norm
\begin{align*}
\|x\|_\gamma:=\|(-A)^\frac\gamma2 x\|,\quad x\in \dot H^\gamma.
\end{align*}
In particular, $\dot H^0=H$.
We will need the following properties of the analytic $\CC_0$-semigroup $S(\cdot)$ (see, e.g., \cite[Theorem 6.13 in Chapter 2]{Paz83}):
\begin{align}\label{ana}
\begin{split}
\|(-A)^\mu S(t)\|_{\LL(H)} \le C t^{-\mu}, \quad
\|(-A)^{-\rho} (S(t)-{\rm Id}_H) \|_{\LL(H)} \le Ct^\rho,
\end{split}
\end{align}
for any $t\in (0,T]$, $\mu\ge 0$ and $\rho \in [0,1]$, where ${\rm Id}_H$ denotes the identity operator in $H$ and $(\LL(H),\|\cdot\|_{\LL(H)})$ denotes the space of bounded linear operators in $H$.

Let $U$ be another separable Hilbert space and ${\bf Q}$ be a self-adjoint, nonnegative definite and bounded linear operator on $U$.
Denote by $U_0:={\bf Q}^\frac12 U$ and $\LL_2^\gamma:=HS(U_0,\dot H^\gamma)$, the Hilbert-Schmidt operator from $U_0$ to $\dot H^\gamma$.
The spaces $H$, $U$ and $\LL_2^\gamma$ are equipped with Borel $\sigma$-algebras $\BB(H)$, $\BB(U)$ and $\BB(\LL_2^\gamma)$, respectively.
Let $W:=\{W(t):\ t\in [0,T]\}$ be a $U$-valued ${\bf Q}$-Wiener process in a stochastic basis $(\Omega,\FFF,(\FFF_t)_{t\in [0,T]},\pp)$, i.e., there exists an eigensystem $\{(q_n, h_n)\}_{n=1}^\infty$ of $\bf Q$ where $\{h_n\}_{n=1}^\infty $ forms an orthonormal basis of $U$ and a sequence of mutually independent Brownian motions $\{\beta_k\}_{n=1}^\infty $ such that (see \cite[Chapter 4]{DZ14})
\begin{align} \label{W}
W(t)=\sum_{n=1}^\infty {\bf Q}^{\frac12} h_n\beta_k(t)
=\sum_{n=1}^\infty \sqrt{q_n} h_n\beta_k(t),\quad t\in [0,T].
\end{align}

\begin{df}\label{df-mild}
A predictable stochastic process $X:[0,T]\times \Omega\rightarrow H$  is called a mild solution of Eq. \eqref{spde} if $X\in L^\infty(0,T;H)$ a.s
and for all $t\in [0,T]$ it holds a.s. that 
\begin{align}\label{mild}
X(t)=S(t)X_0+S*F(X)(t)+S\diamond G(X)(t),
\end{align}
where $S*F(X)$ and $S\diamond G(X)$ denote the deterministic and stochastic convolutions, respectively:
\begin{align*}
&S*F(X)(\cdot):=\int_0^\cdot S(\cdot-r) F(X(r)){\rm d}r, \\
&S\diamond G(X)(\cdot):=\int_0^\cdot S(\cdot-r) G(X(r)){\rm d}W(r).
\end{align*}
We say that $X$ is the unique mild solution of Eq. \eqref{spde} if $Y$ is another solution, then $X$ and $Y$ are stochastically equivalent, i.e.,
$\pp\{X(t)=Y(t)\}=1$, $t\in [0,T]$. 
\end{df}

Let $\theta\ge 0$. 
We use $L^p(\Omega;\CC([0,T]; \dot H^\theta))$ to denote the Banach space consisting of $\dot H^\theta$-valued a.s. continuous stochastic processes $X=\{X(t):\ t\in [0,T]\}$ such that 
\begin{align*}
\|X\|_{L^p(\Omega;\CC([0,T]; \dot H^\theta))}
:=\Big(\ee\Big[\sup_{t\in [0,T]} \|X(t)\|_\theta^p \Big]\Big)^\frac1p
<\infty,
\end{align*}
and 
$L^p(\Omega;\CC^\delta([0,T]; \dot H^\theta))$ with $\delta\in (0,1]$ to denote $\dot H^\theta$-valued a.s. continuous stochastic processes $X=\{X(t):\ t\in [0,T]\}$ such that 
\begin{align*}
\|X\|_{L^p(\Omega;\CC^\delta([0,T]; \dot H^\theta))}
: & =\Big(\ee\Big[\sup_{t\in [0,T]} \|X(t)\|_\theta^p\Big]\Big)^\frac1p \\
&\quad +\bigg(\ee\bigg[\bigg(\sup_{t,s\in [0,T], t\neq s}\frac{\|X(t)-X(s)\|_\theta}{|t-s|^\delta}\bigg)^p\bigg]\bigg)^\frac1p<\infty.
\end{align*}
Our main aim is to find the optimal constants $\delta$ and $\theta$ such that the solution of Eq. \eqref{spde} is in  
$L^p(\Omega;\CC^\delta([0,T];\dot H^\theta))$.
For convenience, throughout $C$ is a generic constant which may be different in each appearance.

\subsection{Main Idea}
\label{sec2.1}

To study the well-posedness and spatial regularity for the solution $X$ of Eq. \eqref{spde}, the main idea of our approach is to use a Burkholder--Davis--Gundy inequality and a weak assumption on the diffusion coefficient $G$ (see Assumption \ref{a2}) to bound the stochastic convolution (see Section \ref{sec3} for more details):
\begin{align*}
\|S\diamond G(X)(t) \|_{L^p(\Omega, \dot H^\theta)}
&\le C \bigg( \int_0^t \|S(t-r) G(X(r)) \|
_{L^p(\Omega, \LL_2^\theta)}^2 {\rm d} r\bigg)^\frac12 \\
&\le C \bigg(\int_0^t K_G^2(t-r) 
\big(1+\|X(r)\|_{L^p(\Omega; \dot H^\theta)} \big)^2 
{\rm d}r\bigg)^\frac12
\end{align*}
for any spatial regularity index $\theta\ge 0$.
Similar argument is applied to the deterministic convolution $S*F(X)$.
Then by H\"older inequality, to bound $\|X(t)\|_{L^p(\Omega; \dot H^\theta)}$ reduces to solve the following type of integral inequality with convolution:
\begin{align}\label{gro}
0\le f(t)\le m(t)+\int_0^t K(t-r) f(r) {\rm d}r,\quad t\in [0,T],
\end{align}
where $f(\cdot)$ is bounded, $m(\cdot)$ is non-decreasing and $K(\cdot)$ is nonnegative and integrable (which may has some singularity at 0) on $[0,T]$.
To overcome this difficulty, we establish a new version of Gr\"onwall inequality with singular kernel, i.e., there exists a constant $\lambda_0$ such that 
$f(t)\le 2e^{\lambda_0 t}m(t)$ (see Lemma \ref{lm-gro}).
Then we obtain the uniform moments' estimation for the solution of Eq. \eqref{spde} under $\|\cdot\|_\theta$-norm (see \eqref{mom} and \eqref{mom-gamma}, respectively).

Using the fixed point argument, a general Lipschitz continuity assumption (see Assumption \ref{a1}) is used to establish the well-posedness as well as the optimal spatial regularity for the solution of Eq. \eqref{spde} (see Section \ref{sec4} for more details).
In this procedure, another difficulty arises from the fact that 
$(\HHH_\theta^p,\|\cdot\|_{\HHH^p})$ (see \eqref{hp} and \eqref{hthetap} for definitions of these two norms) for $\theta>0$ is not a Banach space, while we only assume that the coefficients are Lipschitz continuous in $\|\cdot\|$-norm rather than $\|\cdot\|_\theta$-norm.
This difficulty is a key problem of regularity analysis for semilinear stochastic partial differential equations (SPDEs) and has been pointed out in \cite{JR12(JDE)} and  \cite{Zha07(JFA)}.
To overcome this difficulty, we first utilize the fact that 
$\HHH_\theta^p(M):=\{Z\in \HHH_\theta^p:\ 
\|Z\|_{\HHH_\theta^p}\le M\}$ with $\|\cdot\|_{\HHH^p}$-norm forms a complete metric space for any $M>0$ and $p>1$ (see Lemma \ref{lm-fix}), which
allows us to apply the Banach fixed point theorem to conclude the existence of a unique local solution of Eq. \eqref{spde}.
Then we obtain the global existence by the aforementioned,  uniform a priori estimation.

Our main idea to deal with the trajectory regularity for the solution $X$ of Eq. \eqref{spde} is the factorization formula
\begin{align}\label{wa}
S\diamond G(X)(t)
=\frac{\sin(\pi \alpha)}{\pi}\int_0^t (t-r)^{\alpha-1} S(t-r) G_\alpha(r) {\rm d}r,
\end{align}
where $\alpha\in (0,1)$ and
\begin{align}\label{ga}
G_\alpha(t)
:=\int_0^t (t-r)^{-\alpha} S(t-r) G(X(r)){\rm d}W(r), \quad t\in [0,T].
\end{align}
Similar factorization formula holds for the deterministic convolution 
$S*F(X)$.
To derive the H\"older continuity for the solution of Eq. \eqref{spde}, we give a generalized characterization (see Proposition \ref{prop-hol}) of temporal H\"older continuity of the linear operator $R_\alpha$ defined by 
\begin{align}\label{ra}
R_\alpha f(t):=\int_0^t (t-r)^{\alpha-1} S(t-r) f(r) {\rm d} r,\quad t\in (0,T].
\end{align}
As a consequence of this characterization, we prove the optimal regularity of the Ornstein--Uhlenbeck process $W_A$ (see Corollary \ref{cor-OU}), which generalizes \eqref{da} to the case $\gamma\ge 2\alpha$.
An interesting consequence of the above characterization formulas is that 
we can obtain stronger moments' estimations \eqref{mom-sup} and \eqref{mom-gamma}, which is not a trivial property for the mild solution of Eq. \eqref{spde} under weak assumptions on its data.

\subsection{Main Results}
\label{sec2.2}

To perform our main results, we give the following assumptions on the coefficients $F$ and $G$.

The first assumption is the following Lipschitz-type continuity and linear growth condition, which is the main condition to yield the well-posedness of Eq. \eqref{spde}.

\begin{ap} \label{a1}
There exist two nonnegative, Borel measurable functions $K_F$ and $K_G$ on $[0,T]$ with
\begin{align*}
K_F^0:=\int_0^T K_F(t) {\rm d}t<\infty \quad \text{and}\quad
K_G^0:=\bigg(\int_0^T K_G^2(t) {\rm d}t\bigg)^\frac12 <\infty,
\end{align*}
such that for any $x,y\in H$ and almost every (a.e.) $t\in [0,T]$ it holds that
\begin{align*}
\|S(t) F(x)\| \le K_F(t)(1+\|x\|),  \quad
& \|S(t) (F(x)-F(y))\| \le K_F(t)\|x-y\|, \\
\|S(t) G(x)\|_{\LL_2^0} \le K_G(t)(1+\|x\|),  \quad
& \|S(t) (G(x)-G(y))\|_{\LL_2^0} \le K_G(t)\|x-y\|.
\end{align*}
\end{ap}

To study the spatial regularity for the solution of Eq. \eqref{spde}, we need more growth conditions on $F$ and $G$.
Throughout $\gamma$ is a nonnegative number, which partially characterizes the spatial regularity for the solution of Eq. \eqref{spde}. 

\begin{ap} \label{a2}
There exist nonnegative, Borel measurable functions $K_{F,\gamma}$ and  $K_{G,\gamma}$ on $[0,T]$ with
\begin{align*}
K_F^\gamma:=\int_0^T K_{F,\gamma}(t) {\rm d}t<\infty 
\quad \text{and}\quad
K_G^\gamma:=\bigg(\int_0^T K_{G,\gamma}^2(t) {\rm d}t\bigg)^\frac12<\infty,
\end{align*}
such that for any $z\in \dot H^\gamma$ and a.e. $t\in [0,T]$ it holds that
\begin{align*}
\|S(t) F(z)\|_\gamma
\le K_{F,\gamma}(t)(1+\|z\|_\gamma), \quad
\|S(t) G(z)\|_{\LL_2^\gamma} 
\le K_{G,\gamma}(t)(1+\|z\|_\gamma).
\end{align*}
\end{ap}

In particular, when $\gamma=0$ we set $K_{F,0}=K_F$ and $K_{G,0}=K_G$.

To obtain the temporal regularity for the solution of Eq. \eqref{spde}, we perform the final assumption. 

\begin{ap} \label{a3}
There exists a constant $\alpha\in (1/p, 1/2)$ with $p>2$ such that
\begin{align*}
K_F^{\gamma,\alpha}:=\int_0^T t^{-\alpha} K_{F,\gamma} (t) {\rm d}t<\infty,\quad
K_G^{\gamma,\alpha}:=\bigg(\int_0^T t^{-2\alpha} K_{G,\gamma}^2 (t) {\rm d}t \bigg)^\frac12 <\infty.
\end{align*}
\end{ap}

\begin{rk}\label{rk-com}
Assumptions \ref{a1}--\ref{a3} are weaker than those of \cite[Section 2]{JR12(JDE)} where the authors assumed that $F:H\rightarrow H$ and $G: H\rightarrow \LL_2^0$ are Lipschitz continuous and for some $\gamma\in [0,1)$, $\|G(z)\|_{\LL_2^\gamma}\le C(1+\|z\|_\gamma)$ for any $z\in \dot H^\gamma$.
Indeed, for any $t\in (0,T]$, $\gamma\in [0,1)$ and $z\in \dot H^\gamma$,
\begin{align*}
\|S(t) F(z)\|_\gamma
\le \|(-A)^\frac12 S(t)\|_{\LL(H)} \cdot \|F(z)\|_{\gamma-1} 
\le C\|(-A)^\frac12 S(t)\|_{\LL(H)} (1+\|z\|_\gamma),
\end{align*}
and 
\begin{align*}
\|S(t) G(z)\|_{\LL_2^\gamma}
\le \|S(t)\|_{\LL(H)} \cdot \|G(z)\|_{\LL_2^\gamma}
&\le C\|S(t)\|_{\LL(H)} (1+\|z\|_\gamma).
\end{align*}
Similarly, for any $t\in [0,T]$ and $x,y\in H$ there holds that 
\begin{align*}
&\|S(t) (F(x)-F(y))\| 
\le C\|(-A)^\frac12 S(t)\|_{\LL(H)} \|x-y\|, \\
&\|S(t) (G(x)-G(y))\|_{\LL_2^0}
\le C\|S(t)\|_{\LL(H)} \|x-y\|.
\end{align*}
Set $K_F=K_{F,\gamma}=C\|(-A)^\frac12 S(\cdot)\|_{\LL(H)}$ and
$K_G=K_{G,\gamma}=C\|S(\cdot)\|_{\LL(H)}$.
By the smooth estimation \eqref{ana}, $K_F,K_{F,\gamma}$ are integrable and $K_G(t),K_{G,\gamma}$ are square integrable on $[0,T]$,
which shows Assumptions \ref{a1}-\ref{a2}.
One can also derive Assumption \ref{a3} with $\alpha<1/2$, since 
\begin{align*}
\int_0^T r^{-\alpha} K_{F,\gamma}(t) {\rm d}r
+\int_0^T r^{-2\alpha} K_{G,\gamma}^2(t) {\rm d}r 
\le C \int_0^T \Big(r^{-(\alpha+\frac12)}+r^{-2\alpha}\Big) {\rm d}r<\infty.
\end{align*}
\end{rk}

Our first main result is the following well-posedness result of Eq. \eqref{spde}.

\begin{tm} \label{main1}
Let $p\ge 2$ and $X_0:\Omega\rightarrow H$ be 
$\FFF_0/\BB(H)$-measurable such that $X_0\in L^p(\Omega;H)$. 
Assume that $S(\cdot)$ is a $C_0$-semigroup and Assumptions \ref{a1} holds.
Then Eq. \eqref{spde} possesses a unique mild solution 
$X=\{X(t):\ t\in [0,T]\}$ such that the following statements hold.
\begin{enumerate}
\item [{\rm (1)}]
There exists a constant $C=C(T,p,K_F^0,K_G^0)$ such that
\begin{align}\label{mom}
\sup_{t\in [0,T]}\ee\Big[\|X(t)\|^p\Big]
&\le C\Big(1+\ee\Big[\|X_0\|^p \Big]\Big).
\end{align}

\item [{\rm (2)}]
The solution $X$ is continuous with respect to $\|\cdot\|_{L^p(\Omega;H)}$:
\begin{align}\label{mean}
\lim_{t_1\rightarrow t_2}\ee\Big[\|X(t_1)-X(t_2)\|^p\Big]=0,
\quad t_1,t_2\in [0,T].
\end{align}
\end{enumerate}
\end{tm}

\begin{rk}
To the best of our knowledge, Theorem \ref{main1} is even new for related deterministic PDEs, i.e., Eq. \eqref{spde} with $G=0$, under the minimum Assumption \ref{a1} on $F$.
\end{rk}

Under the conditions of Theorem \ref{main1}, similarly to the additive case as in \cite{DZ14}, one can say nothing about the continuity of the trajectory for the solution $X$ of Eq. \eqref{spde}.
However, if Assumption \ref{a3} holds for $\gamma=0$, we can show that $X$ possesses a continuous version in $H$ by the factorization method even in the case of $C_0$-semigroup.
Moreover, we derive more stronger moments' estimation than \eqref{mom}.

\begin{tm} \label{main2}
In addition to the conditions of Theorem \ref{main1} with $p>2$, assume that Assumption \ref{a3} holds for $\gamma=0$.
Then the mild solution $X$ of Eq. \eqref{spde} belongs to 
$L^p(\Omega;\CC([0,T]; H))$.
Moreover, there exists a constant $C=C(T,p,\alpha,K_F^{0,\alpha}, K_G^{0,\alpha})$ such that 
\begin{align}\label{mom-sup}
\ee\Big[\sup_{t\in [0,T]}\|X(t)\|^p\Big]
&\le C\Big(1+\ee\Big[\|X_0\|^p \Big]\Big).
\end{align}
\end{tm}

Our next main result is the following optimal spatial regularity for the solution of Eq. \eqref{spde}.

\begin{tm} \label{main3}
Let $\gamma>0$, $p\ge 2$ and $X_0:\Omega\rightarrow \dot H^\gamma$ be $\FFF_0/\BB(\dot H^\gamma)$-measurable such that $X_0\in L^p(\Omega;\dot H^\gamma)$. 
Assume that $S(\cdot)$ is an analytic $\CC_0$-semigroup and Assumptions \ref{a1}-\ref{a2} hold.
Then the mild solution $X$ of Eq. \eqref{spde} satisfies the following statements.
\begin{enumerate}
\item [{\rm (1)}]
There exists a constant $C=C(T,p,K_F^\gamma,K_G^\gamma)$ such that
\begin{align}\label{mom-gamma}
\sup_{t\in [0,T]}\ee\Big[\|X(t)\|_\gamma^p\Big]
&\le C\Big(1+\ee\Big[\|X_0\|_\gamma^p \Big]\Big).
\end{align}

\item [{\rm (2)}]
The solution $X$ is continuous with respect to $\|\cdot\|_{L^p(\Omega;\dot H^{\gamma)}}$:
\begin{align}\label{mean-gamma}
\lim_{t_1\rightarrow t_2}\ee\Big[\|X(t_1)-X(t_2)\|_\gamma^p\Big]=0,
\quad t_1,t_2\in [0,T].
\end{align}
\end{enumerate}
\end{tm}

Analogously to Theorem \ref{main2}, we can obtain more stronger moments' estimation than \eqref{mom-sup} and show the a.s. continuity for the solution of Eq. \eqref{spde} in $\dot H^{\gamma}$, under the additional Assumption \ref{a3}.
Moreover, our last main result derives the following optimal trajectory regularity for the solution of Eq. \eqref{spde}.

\begin{tm} \label{main4}
In addition to the conditions of Theorem \ref{main3} with $p>2$, assume that $X_0:\Omega\rightarrow \dot H^\beta$ is $\FFF_0/\BB(\dot H^\beta)$-measurable such that $X_0\in L^p(\Omega;\dot H^\beta)$ and Assumption \ref{a3} holds with $\beta\ge \gamma\ge 0$.
Then the following statements hold.
\begin{enumerate}
\item [{\rm (1)}]
When $\gamma=0$, for any $\delta\in [0,\alpha-1/p)$, $\theta_1\in (0,2\alpha-2/p)$ and $\theta_2\le \beta$ there holds that  
\begin{align}\label{con-0}
X \in
L^p(\Omega;\CC^\delta([0,T]; H) 
& \cup L^p(\Omega;\CC^{\alpha-\frac1p-\frac{\theta_1} 2}([0,T]; \dot H^{\theta_1})) \nonumber \\
& \cap L^p(\Omega;\CC^{\frac{\beta-\theta_2}2\wedge 1}([0,T]; \dot H^{\theta_2})).
\end{align} 

\item [{\rm (2)}]
When $\gamma>0$, for any $\delta\in [0,\alpha-1/p)$, $\theta\in (0,\gamma)$, $\theta_1\in (\gamma, \gamma+2\alpha-2/p)$ and
$\theta_2\le \beta$ there holds that 
\begin{align}\label{con-gamma}
X & \in L^p(\Omega;\CC^\delta ([0,T]; \dot H^{\gamma}) 
\cup  L^p(\Omega;\CC^{\alpha-\frac1p} ([0,T]; \dot H^\theta)) \nonumber \\
& \cup L^p(\Omega;\CC^{\alpha-\frac1p+\frac{\gamma-\theta_1}2} ([0,T]; \dot H^{\theta_1}))
\cap L^p(\Omega;\CC^{\frac{\beta-\theta_2}2\wedge 1}([0,T]; \dot H^{\theta_2})).
\end{align} 
\end{enumerate}
\end{tm}

\begin{rk}
In the analytic $\CC_0$-semigroup case, \eqref{con-0} strengthens the continuity results in Theorem \ref{main2}.
Moreover, \eqref{con-0} and \eqref{con-gamma} show the a.s. continuity for the solution of Eq. \eqref{spde} in $\dot H^\beta$ for $\beta<\gamma+2\alpha-2/p$.
When $\beta\ge \gamma+2\alpha-2/p$, one could not expect that the solution of Eq. \eqref{spde} is a.s. continuous in 
$\dot H^\beta$ due to the optimal regularity of the Ornstein--Uhlenbeck process; see Corollary \ref{cor-OU}.
We also note that \eqref{con-0} and \eqref{con-gamma} show the H\"older regularity for the solution of Eq. \eqref{spde} in $\dot H^\theta$ when 
$\beta\le 2/p$ with $\theta\in [0, \gamma+2\alpha-2/p)\cap [0,\beta]$ 
or $\beta>2/p$ with $\theta\le \beta<\gamma+2\alpha-2/p$.
\end{rk}

\begin{rk}
Theorems \ref{main1}--\ref{main4} establish the well-posedness and optimal trajectory regularity of the solution of Eq. \eqref{spde} for general $\beta$ and $\gamma$ under more general Assumptions \ref{a1}--\ref{a3}, and thus give an answer to Main Problem \ref{MP}.
\end{rk}

Applying our main results in Theorems \ref{main1}--\ref{main4}, we have the following well-posedness and regularity results for Eq. \eqref{spde} under the type of assumptions in \cite{JR12(JDE)}. 

\begin{cor} \label{main-cor}
Let $\beta\ge \gamma\ge 0$, $p\ge 2$ and $X_0:\Omega\rightarrow \dot H^\beta$ be $\FFF_0/\BB(\dot H^\beta)$-measurable such that $X_0\in L^p(\Omega;\dot H^\beta)$. 
Assume that $S(\cdot)$ is a $C_0$-semigroup and $F:H\rightarrow \dot H^{-1}$, $G: H\rightarrow \LL_2^0$ are Lipschitz continuous.

\begin{enumerate}
\item
Eq. \eqref{spde} possesses a unique mild solution $X=\{X(t):\ t\in [0,T]\}$ which belongs to $L^p(\Omega; L^\infty(0,T; H))$ such that \eqref{mean} and \eqref{mom-sup} hold.
If $p>2$, then $X\in L^p(\Omega; \CC([0,T];H))$.
Assume in addition that $S(\cdot)$ is analytic, then 
\begin{align*}
X \in L^p(\Omega; \CC^{\delta_1}([0,T]; \dot H^{\theta_1})) 
\cap L^p(\Omega; \CC^{\frac{\beta-\theta_2}2\wedge 1}([0,T]; \dot H^{\theta_2}))
\end{align*} 
for any $\delta_1,\theta_1,\theta_2\ge 0$ with 
$\delta_1+\theta_1/2<1/2-1/p$ and
$\theta_2\le \beta$.

\item
If $S(\cdot)$ is analytic and $\|F(x)\|_{\gamma-1} \le C(1+\|x\|_\gamma)$, 
$\|G(x)\|_{\LL_2^\gamma}\le C(1+\|x\|_\gamma)$.
Then $X\in L^p(\Omega; L^\infty(0,T; 
\dot H^{\gamma}))$ such that \eqref{mean-gamma} holds.
If $p>2$, then
\begin{align*}
X \in L^p(\Omega; \CC^{\delta_1}([0,T]; \dot H^{\theta_1})) 
\cap L^p(\Omega; \CC^{\frac{\beta-\theta_2}2\wedge 1}([0,T]; \dot H^{\theta_2}))
\end{align*} 
for any $\delta_1,\theta_1,\theta_2\ge 0$ with 
$\delta_1<[1\wedge (\gamma+1-\theta_1)]/2-1/p$
and $\theta_2\le \beta$.
\end{enumerate}
\end{cor}

{\it Proof.}
Taking into account Remark \ref{rk-com}, we note that Assumptions \ref{a1}--\ref{a3} hold with $\alpha<1/2$.
Thus we conclude the first claim by applying Theorems \ref{main1}, \ref{main2} and \ref{main4} and another claim by applying Theorems \ref{main3} and \ref{main4}.
\quad \\

\subsection{Examples}

The main aim of this part is to give several concrete examples which satisfy our main Assumptions \ref{a1}-\ref{a3}.
Our main model is the following second order parabolic SPDE: 
\begin{align}\label{she}\tag{SHE}
\begin{split}
&{\rm d} X(t,\xi)
=(\Delta X(t,\xi)+\nabla \cdot f(X(t,\xi))) {\rm d}t +g(X(t,\xi)) {\rm d}W(t,\xi),\\
&X(t,\xi)=0, \quad (t,\xi)\in [0,T]\times \partial \OOO,\\
&X(0,\xi)=X_0(\xi), \quad \xi\in \OOO,
\end{split}
\end{align}where $\OOO\subset \rr^d$ is a bounded open set with regular boundary. 
Without loss of generality, we assume that $X_0$ is a deterministic function which vanishes on the boundary $\partial \OOO$.

Set $U=H=L^2(\OOO)$ and $A=\Delta$ with domain 
$\text{Dom}(A)=H_0^1(\OOO)\cap H^2(\OOO)$.
Then there exists an eigensystem $\{(\lambda_n,e_n)\}_{n=1}^\infty$ of $-A$: $-A e_n=\lambda_n e_n$, $k\in \nn_+$, where 
$\{\lambda_n\}_{n=1}^\infty$ is in an increasing order and $\{e_n\}_{n=1}^\infty$ forms an orthonormal basis of $H$.
Assume that $f,g:\rr\rightarrow \rr$ are Lipschitz continuous functions with Lipschitz constant $L_f,L_g>0$, i.e., for any $\xi_1$ and $\xi_2\in \rr$ there holds that 
\begin{align} \label{fg}
|f(\xi_1)-f(\xi_2)|\le L_f |\xi_1-\xi_2|,
\quad |g(\xi_1)-g(\xi_2)|\le L_g |\xi_1-\xi_2|.
\end{align}
Let $\{(q_n, h_n)\}_{n=1}^\infty$ be an eigensystem of $\bf Q$ where 
$\{h_n\}_{n=1}^\infty $ forms an orthonormal basis of $H$, and 
$W=\{W(t):\ t\in [0,T]\}$ be an $H$-valued $\bf Q$-Wiener process given by \eqref{W}.
Define the Nemytskii operators $F:H\rightarrow \dot H^{-1}(\OOO)$ and $G:H\rightarrow \LL(H)$, respectively,  by 
\begin{align} \label{FG}
F(x)(\xi):=\nabla \cdot f(x(\xi)),\quad 
G(x) h_n(\xi):=\sqrt{q_n} g(x(\xi)) h_n(\xi),
\end{align}
for $x\in H$, $k\in \nn_+$ and $\xi\in \OOO$.
Then Eq. \eqref{she}  is equivalent to Eq. \eqref{spde} with $F$ and $G$ given by \eqref{FG}.

In the following we will use $(L^\infty(\OOO), \|\cdot\|_{L^\infty(\OOO)})$ to denote the essentially bounded function space and $(\CC^\epsilon(\OOO), \|\cdot\|_{\CC^\epsilon(\OOO)})$ for some $\epsilon\in (0,1)$ to denote the H\"older function space over $\OOO$.

\subsubsection{White Noise}

We begin with the case of white noise.
Assume that $W=\{W(t):\ t\in [0,T]\}$ is an $H$-valued cylindrical Wiener process, i.e., ${\bf Q}={\rm Id}_H$ or equivalently, $q_n=1$ for each $k\in \nn_+$ in \eqref{W}.
In this case, it is known that $G$ defined by \eqref{FG} is not a Lipschitz continuous operator from $H$ to $\LL_2^0$; indeed, $G(H)\nsubseteq \LL_2^0$.
However, we can verify that $F$ and $G$ satisfies Assumptions \ref{a1}--\ref{a3} with $\gamma=0$ (Assumption \ref{a2} reduces to Assumption \ref{a1} when $\gamma=0$).

Let $t\in (0,T]$ and $x,y,z\in H$.
By the definition of $\LL_2^0$-norm and the estimate $\|e_n\|_{L^\infty(\OOO)}\le C \lambda_n^{(d-1)/2}$, $n \in \nn_+$, for the eigensystem of Dirichlet Laplacian (see, e.g. \cite{Gri02(CPDE)}), we get 
\begin{align*}
\|S(t) G(z)\|_{\LL_2^0}^2 
% &=\sum_{n=1}^\infty \|S(t) (g(z) e_n) \|^2
=\sum_{n=1}^\infty e^{-2\lambda_n t} \|g(z) e_n\|^2 
% &\le \sum_{n=1}^\infty e^{-2\lambda_n t} \|e_n\|_{L^\infty(\OOO)}^2 \|g(z) \|^2 \\
\le C \sum_{n=1}^\infty \bigg[\lambda_n^{d-1} e^{-2\lambda_n t} \bigg] (1+\|z\|^2).
\end{align*}
Similarly,
\begin{align*}
\|S(t) (G(x)-G(y)) \|_{\LL_2^0}^2 
\le C \sum_{n=1}^\infty \bigg[\lambda_n^{d-1} e^{-2\lambda_n t} \bigg] \|x-y\|^2.
\end{align*}
Define 
\begin{align} \label{KG}
K_G(t):=C \bigg(\sum_{n=1}^\infty \bigg[\lambda_n^{d-1} e^{-2\lambda_n t} \bigg] \bigg)^\frac12,
\quad t \in (0,T].
\end{align}
By Weyl's law that $\lambda_n\simeq m^{2/d}$ (here $M\simeq N$ means $C_1 N\le M\le C_2 N$ for two nonnegative numbers $C_1$ and $C_2$), we obtain
\begin{align*}
\int_0^T K_G^2(t)  {\rm d}t
\simeq \sum_{n=1}^\infty m^\frac{2(d-2)}d,
\end{align*}
which converges if and only if $d<4/3$.
Thus only for $d=1$, $K_G$ defined by \eqref{KG} is square integrable on $[0,T]$.
Meanwhile, for $\alpha<1/4$,
\begin{align*}
\int_0^T t^{-2\alpha}K_G^2(t)  {\rm d}t
\le C \int_0^T t^{-2\alpha-\frac12} {\rm d}t<\infty.
\end{align*}

On the other hand, for the nonlinear drift term, by the definition \eqref{FG} and the Lipschitz condition \eqref{fg} we get  
\begin{align*}
\|S(t) F(z) \| 
&\le \|(-A)^\frac12 S(t)\|_{\LL(H)} \|F(z) \|_{-1}
\le C t^{-\frac12} (1+\|z\|), \\
&\|S(t) (F(x)-F(y)) \| 
\le C t^{-\frac12} \|x-y\|.
\end{align*}
Define 
\begin{align} \label{KF}
K_F(t):=C t^{-\frac12},
\quad t \in (0,T].
\end{align}
Then the function $K_F$ defined by \eqref{KF} is integrable on $[0,T]$ and for $\alpha<1/2$,
\begin{align*}
\int_0^T t^{-\alpha} K_F(t)  {\rm d}t
\le C \int_0^T t^{-\alpha-\frac12} {\rm d}t<\infty.
\end{align*}

Thus we have shown Assumptions \ref{a1}--\ref{a3} with 
$\gamma=0$ and $\alpha\in (0,1/4)$.
As a result of Theorem \ref{main4} with $X_0\in \dot H^{1/2}$, $\gamma=0$ and $\alpha\in (0,1/4)$, Eq. \eqref{she} driven by an $H$-valued cylindrical Wiener process possesses a unique mild solution in $L^p(\Omega;\CC^\delta([0,T];\dot H^\theta))$ for any $p\ge 1$ and $\delta,\theta\ge 0$ with $\delta+\theta/2<1/4$.

\subsubsection{Colored Noises}

Next we give an example in the case of colored noises which satisfies Assumptions \ref{a1}-\ref{a3} for some $\gamma>0$ and generalizes the examples from \cite[Section 4]{JR12(JDE)}.

Let $\gamma\in (0,1)$, $t\in (0,T]$, $x,y\in H$ and $z\in \dot H^\gamma$. 
For $\gamma\in (0,1/2)$, by the Lipschitz condition \eqref{fg} we have 
$f(z)\in \dot H^\gamma$ and  
$\|f(z) \|_\gamma\le C(1+\|z\|_\gamma)$ for any $z\in \dot H^\gamma$.
This inequality holds true for any $\gamma\in (1/2,1)$ provided that 
$f(0)=0$.
Such additional requirement is due to the characterization of 
$\dot H^\gamma$ (see, e.g., \cite[Appendix (A.46)]{DZ14}):
\begin{align} \label{H-W}
\dot H^\gamma=
\begin{cases}
W^{\gamma,2}(\OOO) 
\quad \text{for}\quad \gamma\in (0,1/2), \\
\big\{x\in W^{\gamma,2}(\OOO):\ x|_{\partial \OOO}=0\big\} 
\quad \text{for}\quad  \gamma\in (1/2,1),
\end{cases}
\end{align}
where $W^{\gamma,2}(\OOO)$ is the Sobolev--Slobodeckij space whose norm is defined by
\begin{align*}
\|X\|_{W^{\gamma,2}(\OOO)}
:=\bigg(\|X\|_{L^2(\OOO)}^2 
+\int_{\OOO}\int_{\OOO} 
\frac{|X(\xi)-X(\eta)|^2}{|\xi-\eta|^{d+2\gamma}}
{\rm d}\xi {\rm d}\eta \bigg)^\frac12.
\end{align*}
It follows by dual argument and the Lipschitz condition \eqref{fg} that 
\begin{align*}
\|S(t) F(z) \|_\gamma 
&\le \|(-A)^\frac12 S(t)\|_{\LL(H)} \|f(z) \|_\gamma 
\le K_{F,\gamma}(t) (1+\|z\|_\gamma), \\
\|S(t) (F(x)-F(y)) \| 
&\le \|(-A)^\frac12 S(t)\|_{\LL(H)} \|f(x)-f(y) \| 
\le C K_{F} (t) \|x-y\|.
\end{align*}
Define 
\begin{align} \label{KF+}
K_{F}(t)=K_{F,\gamma}(t):=C t^{-\frac12},
\quad t \in (0,T].
\end{align}
Then the functions $K_{F}$ and $K_{F,\gamma}$ defined by \eqref{KF+} are integrable on $[0,T]$ and for any $\alpha<1/2$,
\begin{align*}
\int_0^T t^{-\alpha} K_{F,\gamma}(t)  {\rm d}t
\le C \int_0^T t^{-\alpha-\frac12} {\rm d}t<\infty.
\end{align*} 

For the diffusion term, we assume that the eigensystem $\{(q_n, h_n)\}_{n\in \nn_+}$ of ${\bf Q}$ satisfies
\begin{align}\label{con-q0}
{\bf Q}_0:=\sum_{n\in \nn_+}  q_n \|h_n\|_{L^\infty(\OOO)}^2
<\infty.
\end{align}
This condition is valid when ${\bf Q}$ is a trace class operator with uniformly bounded eigenfunctions.
We use the uniform boundedness \eqref{ana}, the Lipschitz condition \eqref{fg} and the assumption \eqref{con-q} to derive
\begin{align*}
\|S(t) (G(x)-G(y)) \|^2_{\LL_2^0}
& \le \|S(t)\|^2_{\LL(H)} \sum_{n\in \nn_+} 
\|(G(x)-G(y)) h_n\|^2  \\
&\le C \sum_{n\in \nn_+} q_n 
\|h_n\|_{L^\infty(\OOO)}^2 \|x-y\|^2
\le C {\bf Q}_0 \|x-y\|^2.
\end{align*}
Similarly, 
\begin{align*}
\|S(t) G(z)\|^2_{\LL_2^\gamma}
\le \|(-A)^\frac\gamma2 S(t)\|^2_{\LL(H)} 
\sum_{n\in \nn_+} \|G(z)h_n \|^2 
\le C  {\bf Q}_0 t^{-\gamma}  (1+ \|z\|)^2.
\end{align*}
Define 
\begin{align} \label{KG+1}
K_G(t):=C,\quad 
K_{G,\gamma}(t):=C t^{-\frac\gamma2},
\quad t \in (0,T].
\end{align}
Then the functions $K_{G}, K_{G,\gamma}$ defined by \eqref{KG+1} are square integrable on $[0,T]$ for any $\gamma<1$ and for 
$\alpha<(1-\gamma)/2$,
\begin{align*}
\int_0^T t^{-2\alpha} K_{G,\gamma}^2(t)  {\rm d}t
\le C \int_0^T t^{-(2\alpha+\gamma)} {\rm d}t<\infty.
\end{align*}
Thus we have shown Assumptions \ref{a1}--\ref{a3} for $\alpha,\gamma>0$ such that $\gamma+2\alpha<1$.
Applying Theorem \ref{main4} with $X_0\in \dot H^1$ and 
$\gamma+2\alpha<1$, Eq. \eqref{she} driven by an $H$-valued ${\bf Q}$-Wiener process $W$ given by \eqref{W} such that \eqref{con-q0} holds possesses a unique mild solution in
$L^p(\Omega;\CC^\delta ([0,T];\dot H^\theta))$
for any $p\ge 1$ and $\delta,\theta\ge 0$ such that $2\delta+\theta<1$.

If more smooth and decay properties on the eigensystem $\{(q_n, h_n)\}_{n\in \nn_+}$ of ${\bf Q}$ are imposed, using Theorem \ref{main4} leads to more regularity for the solution.
Assume that there exists a constant $\epsilon\in (0,1]$ such that 
\begin{align}\label{con-q}
{\bf Q}_\epsilon:=\sum_{n\in \nn_+}  q_n \|h_n\|_{\CC^\epsilon(\OOO)}^2
<\infty.
\end{align}
By the uniform boundedness \eqref{ana}, we get 
\begin{align*}
\|S(t) G(z)\|^2_{\LL_2^\gamma}
\le \|S(t)\|^2_{\LL(H)} 
\sum_{n\in \nn_+} q_n \|g(z)h_n \|_\gamma^2
\le C \sum_{n\in \nn_+} q_n \|g(z)h_n \|_\gamma^2.
\end{align*}
It is shown in \cite[(27) in Section 4]{JR12(JDE)} that  
\begin{align*}
\sum_{n\in \nn_+} q_n \|g(z) h_n\|^2_{W^{\gamma,2}(\OOO)}
\le C \sum_{n\in \nn_+}  q_n \|h_n\|_{\CC^\epsilon(\OOO)}^2 
\|g(z)\|_{W^{\gamma,2}(\OOO)}^2,
\quad \forall~\gamma<\epsilon.
\end{align*}
Then we conclude by the Lipschitz condition \eqref{fg}, the assumption \eqref{con-q} and the characterization \eqref{H-W} that 
\begin{align*}
\|S(t) G(z)\|^2_{\LL_2^\gamma}
\le C {\bf Q}_\epsilon(1+\|z\|_\gamma)^2,
\end{align*}
for any $\gamma<1/2 \wedge \epsilon$ and for any $\gamma\in (0,\epsilon) \setminus \{1/2\}$ provided that $g(0)=0$ or $h_n|_{\partial \OOO}=0$ for all $n\in \nn_+$.
Define 
\begin{align} \label{KG+}
K_{G} (t)=K_{G,\gamma}(t)
:=C {\bf Q}_\epsilon,
\quad t \in (0,T].
\end{align}
Then the functions $K_{G}, K_{G,\gamma}$ defined by \eqref{KG+} are square integrable on $[0,T]$ and for $\alpha<1/2$,
\begin{align*}
\int_0^T t^{-2\alpha} K_{G,\gamma}^2(t)  {\rm d}t
\le C \int_0^T t^{-2\alpha} {\rm d}t<\infty.
\end{align*}

Thus we have shown Assumptions \ref{a1}--\ref{a3} with $\alpha\in (0, 1/2)$ and $\gamma\in (0,1/2 \wedge \epsilon)$ or $\gamma\in (0,\epsilon) \setminus \{1/2\}$ provided $f(0)=g(0)=0$.
Applying Theorem \ref{main4} with $X_0\in \dot H^{3/2}$, 
$\gamma\in (0,1/2 \wedge \epsilon)$ and $\alpha\in (0, 1/2)$, Eq. \eqref{she} driven by an $H$-valued ${\bf Q}$-Wiener process $W$ given by \eqref{W} such that \eqref{con-q} holds for some $\epsilon\in (0,1]$ possesses a unique mild solution in
\begin{align*}
L^p(\Omega;\CC^{\delta_1}([0,T];\dot H^\gamma))
\cup L^p(\Omega;\CC^{\delta_2}([0,T];\dot H^\theta))
\end{align*} 
for any $p\ge 1$, $\delta_1\in (0, 1/2)$, $\theta\in (\gamma,1+\gamma)$,  $\delta_2\in (0, (1+\gamma-\theta)/2)$ and $\gamma\in (0,1/2 \wedge \epsilon)$.
Assume furthermore that $X_0\in \dot H^2$, $f(0)=0$ and $g(0)=0$ or $h_n|_{\partial \OOO}=0$ for all $n\in \nn_+$, then this solution belongs to
\begin{align*}
L^p(\Omega;\CC^{\delta_1}([0,T];\dot H^\gamma))
\cup L^p(\Omega;\CC^{\delta_2}([0,T];\dot H^\theta))
\end{align*} 
for any $p\ge 1$, $\delta_1\in (0, 1/2)$, $\theta\in (\gamma,1+\gamma)$,  
$\delta_2\in (0, (1+\gamma-\theta)/2)$ and $\gamma\in (0,\epsilon) \setminus \{1/2\}$.

\section{Well-posedness and Optimal Spatial Regularity}
\label{sec3}

Our main task in this section is to establish the well-posedness and the optimal spatial regularity for the solution of Eq. \eqref{spde}.

We first establish the well-posedness and uniform $p$-moments' estimation \eqref{mom} for the solution $X$ of Eq. \eqref{spde} under Assumption \ref{a1} (see Theorem \ref{tm-well}).
Then we show that $X$ is continuous in $L^p(\Omega;H)$ (see Proposition \ref{prop-con-mean}).
Combining these results and arguments, we give the proofs of Theorems \ref{main1} and \ref{main3} at the end of this section.

\subsection{Well-posedness}
\label{sec3.1}

For $p\ge 2$, denote by $\HHH^p$ the space of all $H$-valued predictable processes $Y$ defined on $[0,T]$ such that 
\begin{align} \label{hp}
\|Y\|_{\HHH^p}
:=\sup_{t\in [0,T]} \Big(\ee\Big[\|X(t)\|^p \Big]\Big)^\frac1p<\infty.
\end{align}
Note that after identifying stochastic processes which are stochastically equivalent, $(\HHH^p,\|\cdot\|_{\HHH^p})$ becomes a Banach space.

To derive the uniform bounds \eqref{mom} and \eqref{mom-sup} for the solution of Eq. \eqref{spde}, we prove a version of Gr\"onwall inequality with singular kernel.

\begin{lm} \label{lm-gro}
Let $m:[0,T]\rightarrow \rr$ be a non-decreasing and bounded function and $K:[0,T]\rightarrow \rr_+$ be a measurable and nonnegative function such that
\begin{align*}
\alpha_T:=\int_0^T K(r) {\rm d}r<\infty.
\end{align*}
Assume that $f$ is nonnegative and bounded on $[0,T]$ such that 
\begin{align*}
f(t)\le m(t)+\int_0^t K(t-r) f(r) {\rm d}r,\quad t\in [0,T].
\end{align*}
Then there exists a constant $\lambda_0=\lambda_0(T, \alpha_T)$ such that 
\begin{align*}
f(t)\le 2e^{\lambda_0 t}m(t).
\end{align*}
\end{lm}

{\it Proof.}
We extend the functions $f, m, K$ to $\widetilde f, \widetilde m, \widetilde K$, respectively, in $\rr$ by setting them to be $0$ outside $[0,T]$.
Then we get
\begin{align*}
\widetilde f(t)\le \widetilde m(t)+\int_0^t \widetilde K(t-r) \widetilde f(r) {\rm d}r.
\end{align*}
Multiplying the above both sides by $e^{-\lambda t}$ with $\lambda\in [0,\infty)$, we obtain
\begin{align*}
e^{-\lambda t} \widetilde f(t)
\le e^{-\lambda t} \widetilde m(t)
+\int_0^t e^{-\lambda (t-r)} \widetilde K(t-r) e^{-\lambda r} \widetilde f(r) {\rm d}r.
\end{align*}
Set $f_\lambda(t)=e^{-\lambda t} \widetilde f(t)$, 
$m_\lambda(t)=e^{-\lambda t} \widetilde m(t)$ and 
$K_\lambda(t)=e^{-\lambda t} K(t)$, $t\in [0,T]$.
Then we have
\begin{align*}
f_\lambda(t)
\le m_\lambda (t)
+\int_0^t K_\lambda(t-r) f_\lambda(r) {\rm d}r.
\end{align*}
Since $\alpha_T(\lambda):=\int_0^T K_\lambda(t) {\rm d}t$ decreases in $[0,\infty)$ and 
\begin{align*}
\lim_{\lambda\rightarrow 0}\alpha_T(\lambda)
=\alpha_T<\infty,\quad
\lim_{\lambda\rightarrow \infty}\alpha_T(\lambda)
=0,
\end{align*}
there exists a $\lambda_0\in (0,\infty)$ such that 
\begin{align*}
\alpha_T(\lambda_0)=\int_0^T K_{\lambda_0}(t) {\rm d}t<\frac12.
\end{align*}
Thus
\begin{align*}
f_{\lambda_0}(t)
& \le m_{\lambda_0}(t)
+\sup_{r\in [0,t]} f_{\lambda_0}(r) \bigg(\int_0^t K_{\lambda_0}(r) {\rm d}r\bigg)  \\
&\le m_{\lambda_0}(t)
+\frac12 \sup_{r\in [0,t]} f_{\lambda_0}(r),
\end{align*}
from which we get
\begin{align*}
\sup_{r\in [0,t]} f_{\lambda_0}(r)\le 2\sup_{r\in [0,t]} m_{\lambda_0}(r).
\end{align*}
Therefore,
\begin{align*}
e^{-\lambda_0 t} f(t)
\le \sup_{r\in [0,t]} e^{-\lambda_0 t} \widetilde f(t) 
\le 2\sup_{r\in [0,t]} e^{-\lambda_0 t} \widetilde m(t)
\le 2m(t).
\end{align*}
Consequently, we have
\begin{align*}
f(t)\le 2e^{\lambda_0 t}m(t).
\end{align*}
This completes the proof.
\quad \\

\begin{rk}
From the proof we can see that the constant 2 can be replaced by any constant larger than 1.
\end{rk}

\begin{tm} \label{tm-well}
Let $p\ge 2$ and $X_0:\Omega\rightarrow H$ is 
$\FFF_0/\BB(H)$-measurable such that $X_0\in L^p(\Omega;H)$. 
Assume that the linear operator $A$ generates a $C_0$-semigroup and Assumption \ref{a1} holds.
Then Eq. \eqref{spde} possesses a unique mild solution $X$ such that \eqref{mom} holds.
\end{tm}

{\it Proof.}
For $X_0\in L^p(\Omega; H)$ and 
$X\in \HHH^p$ define an operator $\MMM$ by
\begin{align}\label{df-M}
\MMM(X)(t)=S(t)X_0+S*F(X)(t)+S\diamond G(X)(t),
\end{align}
where $t\in [0,T]$. 
We first show that $\MMM$ maps $\HHH^p$ to 
$\HHH^p$. 

By Minkovskii inequality, we get
\begin{align*}
\big\|\MMM(X)\big\|_{\HHH^p}
\le \big\|S(t)X_0\big\|_{\HHH^p}
+\big\|S*F(X)(t)\big\|_{\HHH^p}
+\big\|S\diamond G(X)(t)\big\|_{\HHH^p}.
\end{align*}
By the uniform boundedness of the semigroup $S$, we set 
\begin{align} \label{df-Mt}
M(t):=\sup_{r\in [0,t]} \|S(r)\|,\quad t\in [0,T].
\end{align}
Then 
\begin{align*}
\big\|S(\cdot)X_0\big\|_{\HHH^p}
\le M_T  \|X_0\|_{L^p(\Omega;H)}.
\end{align*}
By Minkovskii inequality and Assumption \ref{a1}, we get
\begin{align*}
\| S*F(X) \|_{\HHH^p} 
&\le \sup_{t\in [0,T]} \int_0^t \|S(t-r) F(X(r))\|_{L^p(\Omega;H)} {\rm d}r \\
&\le \sup_{t\in [0,T]} \int_0^t K_F(t-r) (1+\|X(r)\|_{L^p(\Omega;H)}) {\rm d}r \\
&\le \bigg(\int_0^T K_F(r){\rm d} r\bigg) \bigg(1+\|X\|_{\HHH^p}\bigg).
\end{align*}
For the stochastic convolution, applying Burkholder--Davis--Gundy inequality and Assumption \ref{a1}, we obtain
\begin{align*}
\ee\bigg[\|S\diamond G(X)(t) \|^p\bigg] 
&\le \bigg(\int_0^t \|S(t-r) G(X(r)) \|_{L^p(\Omega, \LL_2^0)}^2  {\rm d} r\bigg)^\frac p2 \\
&\le \bigg(\int_0^t K_G^2(t-r) 
\big(1+\|X(r)\|_{L^p(\Omega; H)} \big)^2
{\rm d}r\bigg)^\frac p2 \\
&\le \bigg(\int_0^t K_G^2(r)  {\rm d}r\bigg)^\frac p2
\Big(1+\|X\|_{\HHH^p} \Big)^p.
\end{align*}
Then 
\begin{align*}
\|S\diamond G(X)\|_{\HHH^p} 
\le \bigg(\int_0^T K_G^2(r)  {\rm d}r\bigg)^\frac12
\Big(1+\|X\|_{\HHH^p} \Big).
\end{align*}
Combining the above estimates, we have
\begin{align*}
\big\|\MMM(X)\big\|_{\HHH^p}
\le M_T  \|X_0\|_{L^p(\Omega;H)}
+N_T \big(1+\|X\|_{\HHH^p} \big),
\end{align*}
where $N(t)$ is the non-decreasing, continuous function defined by
\begin{align*}
N(t)=
\int_0^t K_F(r){\rm d} r+\bigg(\int_0^t K_G^2(r)  {\rm d}r\bigg)^\frac12,\quad t\in [0,T].
\end{align*}
Thus $\big\|\MMM(X)\big\|_{\HHH^p}<\infty$ and $\MMM$ maps $\HHH^p$ to $\HHH^p$. 

Next we show that $\MMM$ is a contraction. 
To this end, we introduce the norm
\begin{align}\label{df-hpu}
\|Y\|_{\HHH^{p,u}}
:=\sup_{t\in [0,T]} e^{-ut}\Big(\ee\Big[\|X(t)\|^p \Big]\Big)^\frac1p,
\end{align}
which is equivalent to $\|\cdot\|_{\HHH^p}$ for any $u>0$.
Then for $X^1,X^2\in \HHH^{p,u}$, previous arguments yield that
\begin{align*}
& \|\MMM(X^1)(t)-\MMM(X^2)(t)\|_{L^p(\Omega;H)} \\
&\le \bigg\|\int_0^t S(t-r) (F(X^1(r))-F(X^2(r))) {\rm d}r 
\bigg\|_{L^p(\Omega;H)}\\
&\quad +\bigg\|\int_0^t S(t-r) (G(X^1(r))-G(X^2(r))){\rm d}W(r)\bigg\|_{L^p(\Omega;H)} \\
&\le \int_0^t K_F(t-r) \|X^1(r)-X^2(r)\|_{L^p(\Omega;H)} {\rm d}r \\
&\quad +\bigg(\int_0^t K_G^2(t-r) \|X^1(r)-X^2(r)\|_{L^p(\Omega;H)}^2 {\rm d}r \bigg)^\frac12 \\
&\le \bigg(\int_0^t e^{u r} K_F(t-r){\rm d}r 
+\bigg(\int_0^t e^{2u r} K_G^2(t-r) {\rm d}r\bigg)^\frac12 \bigg) 
\|X^1-X^2\|_{\HHH^{p,u}}.
\end{align*}
Then
\begin{align*}
\|\MMM(X_1)-\MMM(X_2)\|_{\HHH^{p,u}}
\le N_T(u) \|X^1-X^2\|_{\HHH^{p,u}},
\end{align*}
where 
\begin{align*}
N_T(u)
&=\sup_{t\in [0,T]}\Bigg[ e^{-u t}\bigg(\int_0^t e^{u r} K_F(t-r){\rm d}r 
+\bigg(\int_0^t e^{2u r} K_G^2(t-r) {\rm d}r\bigg)^\frac12 \bigg)\Bigg] \\
&=\int_0^T e^{-ur} K_F(r){\rm d}r 
+\bigg(\int_0^T e^{-2u r} K_G^2(r) {\rm d}r\bigg)^\frac12.
\end{align*}
It is clear that the function $N_T:\rr_+\rightarrow \rr_+$ is non-increasing and continuous with $N_T(0)=N_T<\infty$ and $N_T(\infty)=0$.
Thus there exists a sufficiently large $u^*\in \rr_+$ such that $N_T(u^*)<1$.
As a consequence, the operator $\MMM$ is a strict contraction in 
$(\HHH^p,\|\cdot\|_{\HHH^{p,u^*}})$, which shows the existence and uniqueness of a mild solution of Eq. \eqref{spde} such that
\begin{align}\label{mom-0}
\sup_{t\in [0,T]} \ee\Big[\|X(t)\|^p\Big]<\infty.
\end{align}
The existence of a predictable version is a consequence of \cite[Proposition 3.6]{DZ14}.

It remains to prove the estimation \eqref{mom}.
Previous idea implies the following estimation:
\begin{align*}
\|X(t)\|_{L^p(\Omega;H)} 
&\le M(t)  \|X_0\|_{L^p(\Omega;H)}+N(t) \\
&\quad +\int_0^t K_F(t-r) \|X(r)\|_{L^p(\Omega;H)} {\rm d}r  \\
&\quad +\bigg(\int_0^t K_G^2(t-r) \|X(r)\|_{L^p(\Omega; H)}^2 {\rm d}r\bigg)^\frac12.
\end{align*}
Then by H\"older inequality, we have
\begin{align*}
\|X(t)\|_{L^p(\Omega;H)}^2 
&\le 3 \Big(M(t)  \|X_0\|_{L^p(\Omega;H)}
+N(t)\Big)^2  \\
&\quad +3 \bigg(\int_0^t K_F(r) {\rm d}r\bigg)
\bigg(\int_0^t K_F(t-r) \|X(r)\|_{L^p(\Omega;H)}^2 
{\rm d}r\bigg)  \\
&\quad +3 \int_0^t K_G^2(t-r) \|X(r)\|_{L^p(\Omega; H)}^2 {\rm d}r.
\end{align*}
Set for $t\in [0,T]$
\begin{align*}
m(t):=3 \big(M(t)  \|X_0\|_{L^p(\Omega;H)}+N(t)\big)^2,  \quad
K(t):=3  \big(K_F^0 K_F(t)+K_G^2(t) \big).
\end{align*}
It is clear that $m$ is non-decreasing and bounded, $K$ is integrable on $[0,T]$ and  
\begin{align*}
\|X(t)\|_{L^p(\Omega;H)}^2 
&\le m(t)+\int_0^t K(t-r) \|X(r)\|_{L^p(\Omega;H)}^2 {\rm d}r.
\end{align*}
Applying the uniform boundedness \eqref{mom-0} and Lemma \ref{lm-gro}, we conclude \eqref{mom}.
\quad \\

\subsection{$L^p(\Omega)$-Continuity}
\label{sec3.2}

Under the conditions of Theorem \ref{tm-well}, we can show that the solution $X$ of Eq. \eqref{spde} is continuous with respect to $\|\cdot\|_{L^p(\Omega;H)}$.

\begin{prop} \label{prop-con-mean} 
Assume that the assumptions of Theorem \ref{tm-well} hold.
Then for any $t_1,t_2\in [0,T]$ there holds that 
\begin{align}\label{cor-con0}
\lim_{t_1\rightarrow t_2}\ee\Big[\|X(t_1)-X(t_2)\|^p\Big]=0.
\end{align}
\end{prop}

{\it Proof.}
Without loss of generality, assume that $0\le t_1<t_2\le T$.
Due to the strong continuity of the 
$C_0$-semigroup $S(t)$: 
\begin{align}
\label{c0}
(S(t)-{\rm Id}_H) x\rightarrow 0\ \text{in}\ H
\ \text{as}\ t\rightarrow 0, \quad \forall\ x\in H,
\end{align} 
the term $S(\cdot)X_0$ is continuous in $L^p(\Omega;H)$:
\begin{align}\label{con-mean1}
& \lim_{t_1\rightarrow t_2}
\ee\Big[\|S(t_1) X_0-S(t_2) X_0 \|^p\Big] \nonumber \\
&=\lim_{t_1\rightarrow t_2}
\ee\Big[\|(S(t_2-t_1)-{\rm Id}_H) S(t_1)X_0 \|^p\Big]
=0.
\end{align}

Next we consider the stochastic convolution $S\diamond G(X)$.
By H\"older and Burkholder--Davis--Gundy inequalities, we get
\begin{align*}
&\ee\Big[\|S\diamond G(X)(t_1)-S\diamond G(X)(t_2) \|^p\Big] \\
&\le \bigg( \int_0^{t_1} \|(S(t_2-t_1)-{\rm Id}_H) S(t_1-r) G(X(r))\|_{L^p(\Omega;\LL_2^0)}^2 {\rm d} r \bigg)^\frac p2  \\
&\quad +\bigg(\int_{t_1}^{t_2} \|S(t_2-r) G(X(r))\|_{L^p(\Omega;\LL_2^0)}^2 {\rm d} r \bigg)^\frac p2
=: I_1+I_2.
\end{align*}
For the first term, by the uniform boundedness of the $C_0$-semigroup $S(t)$ and the uniformly boundedness \eqref{mom} of $X$, we get 
\begin{align*}
I_1\le C \bigg(\int_0^{t_1} K_G^2(r) {\rm d} r\bigg)^\frac p2
\bigg(1+\|X\|_{\HHH^p}\bigg)^p
<\infty.
\end{align*}
Then $I_1$ tends to 0 as $t_1\rightarrow t_2$ by the strong continuity \eqref{c0} of the $C_0$-semigroup $S(t)$ and Lebesgue dominated convergence theorem.
For the second term, we have
\begin{align*}
I_2\le \bigg(\int_0^{t_2-t_1} K_G^2 (r) {\rm d}r\bigg)^\frac p2
\Big(1+\|X\|_{\HHH^p} \Big)^p
\rightarrow 0\quad \text{as}\quad  t_1\rightarrow t_2
\end{align*}
by Lebesgue dominated convergence theorem.
Therefore,
\begin{align}\label{con-mean2}
\lim_{t_1\rightarrow t_2}
\ee\Big[\|S\diamond G(X)(t_1)-S\diamond G(X)(t_2) \|^p\Big]=0.
\end{align}
Similar arguments can handle the deterministic convolution $S*F(X)$:
\begin{align}\label{con-mean3}
\lim_{t_1\rightarrow t_2}
\ee\Big[\|S*F(X)(t_1)-S*F(X)(t_2) \|^p\Big]=0.
\end{align}
Combining the estimations \eqref{con-mean1}--\eqref{con-mean3}, we derive
\eqref{mean}.
\quad \\

\subsection{Proof of Theorems \ref{main1} and  \ref{main3}}
\label{sec3.3}

In this part, we prove Theorems \ref{main1} and \ref{main3}.

%%%%%%%%%%%%%%%%%%%%%%%%%%%%%%%%%%%%%%%%%

{\it Proof.}[Proof of Theorem \ref{main1}]
Combining Theorem \ref{tm-well} and Proposition \ref{prop-con-mean}, we conclude Theorem \ref{main1}.
\quad \\

%%%%%%%%%%%%%%%%%%%%%%%%%%%%%%%%%%%%%%%%%

To study the spatial regularity for the solution of Eq. \eqref{spde} and prove Theorem \ref{main3}, for $\theta>0$ and $p\ge 2$, we denote by
$\HHH_\theta^p$ the space of all $H$-valued predictable processes $Y$ defined on $[0,T]$ such that 
\begin{align} \label{hthetap}
\|Y\|_{\HHH_\theta^p}
:=\sup_{t\in [0,T]} \Big(\ee\Big[\|X(t)\|_\theta^p \Big]\Big)^\frac1p<\infty.
\end{align}
Unlike the proof of Theorem \ref{tm-well} where we used the fact that 
$(\HHH^p,\|\cdot\|_{\HHH^p})$ is a Banach space, 
$(\HHH_\theta^p,\|\cdot\|_{\HHH^p})$ for $\theta>0$ does not forms a Banach space.
However, the following result shows that 
\begin{align} \label{df-hp}
\HHH_\theta^p(M):=\{Z\in \HHH_\theta^p:\  \|Z\|_{\HHH_\theta^p}\le M\}
\end{align} 
with norm $\|\cdot\|_{\HHH^p}$ is a complete metric space for any $M>0$ and $p\ge 2$.

\begin{lm} \label{lm-fix}
For any $M>0$, $p>1$ and $\theta\ge 0$, the space 
$\HHH_\theta^p(M)$ defined by \eqref{df-hp} with norm $\|\cdot\|_{\HHH^p}$ is a complete metric space.
\end{lm}

{\it Proof.}
Let $M>0$, $p\ge 1$ and $\theta\ge 0$.
Assume that $\{u_n\}_{n\in \nn_+}\subset \HHH_\theta^p(M)$ and $u_n\rightarrow u$ in $\HHH^p$ as $n\rightarrow \infty$.
Then $\{u_n\}_{n\in \nn_+}$ is uniformly bounded in $\HHH_\theta^p$ by $M$ and thus there exists a subsequence, which we still denote by 
$\{u_n\}_{n\in \nn_+}$, such that $u_n(t)\rightarrow u(t)$ in $L^p(\Omega; H)$ for a.e. $t\in [0,T]$.

Since for each $p>1$ and $\theta\ge 0$ the space 
$L^p(\Omega; \dot H^\theta)$ is reflexible and 
\begin{align*}
 \big(L^p(\Omega; \dot H^\theta), \|\cdot \|_{L^p(\Omega; \dot H^\theta)} \big)\hookrightarrow  \big(L^p(\Omega; H), \|\cdot \|_{L^p(\Omega; H)} \big),
\end{align*} 
we conclude by \cite[Theorem 1.2.5]{Caz03} that the limit $u$ belongs to $\HHH_\theta^p$ such that 
\begin{align*}
\|u \|_{L^\infty(0,T; L^p(\Omega; \dot H^\theta))} 
\le \liminf_{n\rightarrow \infty} \|u_n\|_{L^\infty(0,T; L^p(\Omega; \dot H^\theta))}
\le M.
\end{align*} 
This shows that $u\in \HHH_\theta^p(M)$ and completes the proof.
\quad \\

Lemma \ref{lm-fix} allows us to apply the Banach fixed point theorem to conclude the existence of a unique local solution of Eq. \eqref{spde}.
Then we prove the global existence by a uniform a priori estimation.\\

{\it Proof.}[Proof of Theorem \ref{main3}]
Let $X\in \HHH_\gamma^p$ and $X^1, X^2\in \HHH^p$.
Using similar arguments as in the proof of Theorem \ref{tm-well} yields that the operator $\MMM$ defined by \eqref{df-M} satisfies 
\begin{align*}
\big\|\MMM(X)\big\|_{\HHH_\gamma^p}
\le M_T  \|X_0\|_{L^p(\Omega; \dot H^\gamma)}
+ & N_\gamma(T) \big(1+\|X\|_{\HHH_\gamma^p} \big), \\
\|\MMM(X_1)-\MMM(X_2)\|_{\HHH^p}
\le N_\gamma& (T) \|X^1-X^2\|_{\HHH^p},
\end{align*}
where $M(\cdot)$ is defined by \eqref{df-Mt} and $N_\gamma$ is a non-decreasing, continuous function defined by
\begin{align*}
N_\gamma(t):=
\int_0^t K_{F,\gamma}(r){\rm d} r+\bigg(\int_0^t K_{G,\gamma}^2(r)  {\rm d}r\bigg)^\frac12,\quad t\in [0,T].
\end{align*}

Since $N_\gamma$ is non-decreasing and continuous with $N_\gamma(0)=0$, there exists a small enough $T$ such that $N_\gamma(T)<1$.
Taking $M$ sufficiently large such that 
\begin{align*}
M\ge \frac{M_T \|X_0\|_{L^p(\Omega; \dot H^\gamma)}+N_\gamma(T)}{1-N_\gamma(T)},
\end{align*}
we conclude that $\MMM$ maps $\HHH_\gamma^p(M)$ to 
$\HHH_\gamma^p(M)$ and is a contraction under the $\|\cdot\|_{\HHH^p}$-norm for sufficiently small time $T$.
By Lemma \ref{lm-fix} and the Banach fixed point theorem, given any $T>0$ there exists a deterministic time $\tau\in (0,T)$ satisfying $N_\gamma(\tau)<1$ such that Eq. \eqref{spde} possesses a unique local mild solution $\{u(t): t\in [0,\tau]\}$ which possesses a predictable version such that
\begin{align}\label{mom-gamma0}
\sup_{t\in [0,\tau]} \ee\Big[\|X(t)\|_\gamma^p\Big]<\infty.
\end{align}

It remains to prove the uniform a priori estimation \eqref{mom-gamma} to conclude the global existence for the solution of Eq. \eqref{spde}.
Let $t\in [0,\tau]$.
Similar arguments as in the proof of Theorem \ref{tm-well} imply the following estimation:
\begin{align*}
\|X(t)\|_{L^p(\Omega;\dot H^\gamma)} 
&\le M(t)  \|X_0\|_{L^p(\Omega;\dot H^\gamma)}+N_\gamma(t) \\
&\quad +\int_0^t K_{F,\gamma}(t-r) \|X(r)\|_{L^p(\Omega;\dot H^\gamma)} {\rm d}r  \\
&\quad +\bigg(\int_0^t K_{G,\gamma}^2(t-r) \|X(r)\|_{L^p(\Omega; \dot H^\gamma)}^2 {\rm d}r\bigg)^\frac12.
\end{align*}
Then by H\"older inequality, we obtain
\begin{align*}
\|X(t)\|_{L^p(\Omega;\dot H^\gamma)}^2 
&\le m_\gamma(t)+\int_0^t K_\gamma(t-r) \|X(r)\|_{L^p(\Omega;\dot H^\gamma)}^2 {\rm d}r,
\end{align*}
where 
\begin{align*}
m_\gamma(t):=3 \big(M(t)  \|X_0\|_{L^p(\Omega;\dot H^\gamma)}
+N_\gamma(t) \big)^2,  \quad
K_\gamma(t):=3 \big(K^\gamma_F K_{F,\gamma}(t)+K_{G,\gamma}^2(t)\big).
\end{align*}
It is clear from Assumptions \ref{a1}-\ref{a2} that $m_\gamma$ is non-decreasing and bounded and $K_\gamma$ is integrable on $[0,T]$.
Then applying Lemma \ref{lm-gro}, we conclude by the boundedness \eqref{mom-gamma0} that 
there exists a constant $C=C(T,p,K_F^\gamma,K_G^\gamma)$ independent of 
$\tau$ such that the aforementioned local solution satisfies the following a priori estimation:
\begin{align*}
\sup_{t\in [0,\tau]} \ee\Big[\|X(t)\|_\gamma^p\Big] 
\le C\Big(1+\ee\Big[\|X_0\|_\gamma^p \Big]\Big).
\end{align*}
Since the above constant $C$ is independent of $\tau$, Eq. \eqref{spde} exists a unique solution on $[0,T]$ such that \eqref{mom-gamma} holds.

To prove \eqref{mean-gamma}, set $t_1<t_2$ without loss of generality.
Let us note that it follows from the proof of Proposition \ref{prop-con-mean} that 
\begin{align*}
&\ee\Big[\|S\diamond G(X)(t_1)-S\diamond G(X)(t_2)\|_\gamma^p\Big] \\
&\le \bigg( \int_0^{t_1} \bigg(\ee\bigg[\big\|(S(t_2-t_1)-I) S(t_1-r) G(X(r))\big\|_\gamma^p\bigg]\bigg)^\frac 2p {\rm d} r \bigg)^\frac p2 \\
&\quad +\bigg(\int_0^{t_2-t_1} K_{G,\gamma}^2 (r) {\rm d}r\bigg)^\frac p2
\bigg(1+\|X\|_{\HHH^p_\gamma}^p \bigg),
\end{align*}
which tends to 0 as $t_1\rightarrow t_2$ by strong continuity \eqref{c0} of the 
$C_0$-semigroup $S(t)$ and Lebesgue dominated convergence theorem.
Similar arguments can handel the deterministic convolution $S*F(X)$ and the term $S(\cdot)X_0$:
\begin{align*}
&\lim_{t_1\rightarrow t_2}
\ee\Big[\|S\diamond G(X)(t_1)-S\diamond G(X)(t_2) \|_\gamma^p\Big]=0,\\
&\lim_{t_1\rightarrow t_2}
\ee\Big[\|S(t_1) X_0-S(t_2) X_0\|_\gamma^p\Big]=0.
\end{align*}
This complete the proof of \eqref{mean-gamma} and thus the proof of Theorem \ref{main3}.
\quad \\

\section{Optimal Trajectory Regularity}
\label{sec4}

Now we consider the trajectory regularity for the solution of Eq. \eqref{spde} in $\dot H^\theta$ for some $\theta\ge 0$.
The main tool is the factorization method introduced in \cite{DKZ87(STO)}.

To derive more temporal regularity of $X$, we generalize a characterization of   the temporal H\"older Continuity for the linear operator $G_\alpha$ defined by \eqref{ra} in \cite[Proposition 5.14]{DZ14} (see Proposition \ref{prop-hol}).
Then we obtain the optimal temporal regularity of $X$ by this characterization (see Theorem \ref{prop-hol}) and thus prove Theorem \ref{main4}.

\subsection{Proof of Theorem \ref{main2}}
\label{sec4.1}

Factorization method is a powerful tool to show the existence of a continuous version of the solution of an SEE.
It is first introduced by G. Da prato, S. Kwapie\u{n} $\&$ J. Zabczyk \cite{DKZ87(STO)} to the stochastic setting; see also \cite[Proposition 5.9]{DZ14}.

We begin with a continuity characterization of $R_\alpha$ defined by \eqref{ra}.

\begin{lm}\label{lm-con}
Let $S(\cdot)$ be a $C_0$-semigroup generated by $A$.
Assume that $p>1$, $\rho\ge 0$, $\alpha>1/p+\rho$ and $E_1$, $E_2$ are Banach spaces such that 
\begin{align*}
\|S(t) x\|_{E_1}\le C t^{-\rho}\|x\|_{E_2},\quad 
t\in (0,T],\ x\in E_2.
\end{align*}
Then $R_\alpha$ defined by \eqref{ra} is a bounded linear operator from $L^p(0,T; E_2)$ to $\CC([0,T]; E_1)$.
\end{lm}

Now we can prove Theorem \ref{main2} by the above lemma.

{\it Proof.}[Proof of Theorem \ref{main2}]
The property \eqref{c0} and Lebesgue dominated convergence theorem yield that
\begin{align*}
\|S(t_2)X_0-S(t_1)X_0\|
&=\|(S(t_2-t_1)-{\rm Id}_H) S(t_1) X_0\| \rightarrow 0
\end{align*}
as $t_1\rightarrow t_2$.
Thus $S(\cdot) X_0\in \CC([0,T]; H)$.

By Burkholder--Davis--Gundy inequality, we get
\begin{align*}
\| G_\alpha(t) \|_{L^p(\Omega;H)}  
&\le \bigg(\int_0^t \|(t-r)^{-\alpha}  S(t-r) G(X(r)) \|_{L^p(\Omega;\LL_2^0)}^2 {\rm d} r \bigg)^\frac12 \\
&\le \bigg(\int_0^t r^{-2\alpha} K_G^2(r) {\rm d}r  \bigg)^\frac12
\bigg(1+\|X\|_{\HHH^p}\bigg).
\end{align*}
Then by Fubini theorem, we get
\begin{align*}
&\ee \bigg[\| G_\alpha(t) \|_{L^p(0,T;H)}^p  \bigg]  
=\int_0^T \ee\bigg[\|G_\alpha(t)\|^p \bigg]  {\rm d}t  \\
&\le \bigg[\int_0^T \bigg(\int_0^t r^{-2\alpha} K_G^2(r) {\rm d}r  \bigg)^\frac p2 {\rm d}t \bigg] \bigg(1+\|X\|_{\HHH^p}\bigg)^p
<\infty.
\end{align*}
This shows that $G_\alpha\in L^p(0,T; H)$ a.s.
Applying Lemma \ref{lm-con} with $E_1=E_2=H$ and $\rho=0$, we have that 
$S\diamond G(X)\in \CC([0,T]; H)$.
Similar argument yields that $S*F(X)\in \CC([0,T]; H)$.
Combining the continuity of $S(\cdot)X_0$, $S*F(X)$ and $S\diamond G(X)$, we complete the proof of the continuity of $X$ in $H$.

For the term $S(t) X_0$, we have
	\begin{align*}
	\ee\Big[\sup_{t\in [0,T]} \|S(t) X_0\|^p \Big]
	\le M_T^p \ee\Big[\|X_0\|^p \Big].
	\end{align*}
	
	By the factorization formula \eqref{wa} for $S\diamond G(X)$ with $G_\alpha$ being given by \eqref{ga}, we have
	\begin{align*}
	& \sup_{t\in [0,T]} \big\|S\diamond G(X)(t)\big\|^p \\
	&\le \bigg(\frac{\sin(\pi \alpha)}{\pi}\bigg)^p
	\bigg(\int_0^T r^{(\alpha-1)p'} {\rm d}r\bigg)^\frac p{p'} 
	\bigg(\sup_{t\in [0,T]}\int_0^t \|S(t-r) G_\alpha(r)\|^p {\rm d}r\bigg) \\
	&\le T^{\alpha p-1} \bigg(\frac{\sin(\pi \alpha)}{\pi}\bigg)^p
	M_T^p \bigg(\int_0^T \|G_\alpha (t)\|^p {\rm d}t\bigg).
	\end{align*}
	
	On the other hand, by Fubini theorem and Burkerholder--Davis--Gundy inequality, we obtain
	\begin{align*}
	& \ee\bigg[\int_0^T \|G_\alpha(t)\|^p {\rm d}t\bigg] 
	=\int_0^T \ee\bigg[\|G_\alpha(t)\|^p \bigg] {\rm d}t  \\
	&\le \int_0^T \bigg(\int_0^t (t-r)^{-2\alpha} 
	\|S(t-r) G(X(r))\|_{L^p(\Omega;\LL_2^0)}^2
	{\rm d}r\bigg)^\frac p2 {\rm d}t \\
	&\le \bigg(\int_0^T \bigg(\int_0^t r^{-2\alpha} K_G^2(r) {\rm d}r\bigg)^\frac p2{\rm d}t \bigg)
	\bigg(1+\|X\|_{\HHH^p}^2\bigg)^\frac p2 \\
	&\le 2^{\frac p2-1} T \Big(K_G^{\gamma,\alpha} \Big)^p  
	\Big(1+\|X\|_{\HHH^p}^p\Big).
	\end{align*}
	Thus we get
	\begin{align*}
	\ee\bigg[\sup_{t\in [0,T]} \big\|S\diamond G(X)(t)\big\|^p\bigg]
	\le 2^{\frac p2-1}  T^{\alpha p} M_T^p  \Big(K_G^{\gamma,\alpha} \Big)^p  
	\bigg(\frac{\sin(\pi \alpha)}{\pi}\bigg)^p
	\Big(1+\|X\|_{\HHH^p}^p\Big).
	\end{align*}
	
	Similarly, we have
	\begin{align*}
	& \ee\bigg[\sup_{t\in [0,T]} \big\|S*F(X)(t)\big\|^p \bigg] \\
	&\le T^{\alpha p-1} M^p_T \bigg(\frac{\sin(\pi \alpha)}{\pi}\bigg)^p
	\bigg(\int_0^T  \ee\bigg[\bigg\|\int_0^t (t-r)^{-\alpha} S(t-r) F(X(r)){\rm d} r \bigg\|^p\bigg] {\rm d}t\bigg) \\
	&\le T^{\alpha p-1} M^p_T  \bigg(\frac{\sin(\pi \alpha)}{\pi}\bigg)^p
	\bigg(\int_0^T \bigg(\int_0^t r^{-\alpha} K_F(r) {\rm d}r\bigg)^p {\rm d}t \bigg)
	\bigg(1+\|X\|_{\HHH^p}\bigg)^p \\
	&\le 2^{p-1}  T^{\alpha p} M_T^p \Big(K_F^{\gamma,\alpha} \Big)^p
	\bigg(\frac{\sin(\pi \alpha)}{\pi}\bigg)^p
	\Big(1+\|X\|_{\HHH^p}^p\Big).
	\end{align*}
	
	Combining the above estimations, we obtain by H\"older inequality that 
	\begin{align*}
	\ee\Big[\sup_{t\in [0,T]} \|X(t)\|^p \Big]
	\le C \Big(1+\ee\Big[\|X_0\|^p \Big]+\|X\|_{\HHH^p}^p \Big),
	\end{align*}
	from which and \eqref{mom} we conclude \eqref{mom-sup}.
\quad \\

\subsection{H\"older Continuity Criterion}
\label{sec4.2}

To deduce more temporal regularity of the deterministic and stochastic convolutions, one needs to assume that $S(\cdot)$ is an analytic $\CC_0$-semigroup generated by $A$.
From now on we assume that the linear operator $A$ generates an analytic $\CC_0$-semigroup such that \eqref{ana} hold.

We have the following characterization of temporal H\"older continuity of the linear operator $R_\alpha$ defined by \eqref{ra}.
The case $\beta=0$ was derived in \cite[Proposition 5.14]{DZ14}.
We give a self-contained proof for completeness.

\begin{prop} \label{prop-hol}
Let $p>1$, $1/p<\alpha<1$ and $\rho,\theta,\delta\ge 0$.
Then $R_\alpha$ defined by \eqref{ra} 
is a bounded linear operator from $L^p(0,T; \dot H^\rho)$ to $\CC^{\delta}([0,T]; \dot H^\theta)$ when $\alpha,\rho,\theta,\delta$ satisfy one of the following conditions:
\begin{enumerate}
\item
$\delta=\alpha-1/p-(\theta-\rho)/2$ when $\theta>\rho$ and 
$\alpha>(\theta-\rho)/2+1/p$;

\item
$\delta<\alpha-1/p$ when $\theta=\rho$;

\item
$\delta=\alpha-1/p$ when $\theta<\rho$.
\end{enumerate}
\end{prop}

{\it Proof.}
Let $0\le t_1<t_2\le T$ and $f\in L^p([0,T]; \dot H^\rho)$.
Then 
\begin{align*}
& \|R_\alpha f(t_2)-R_\alpha  f(t_1) \|_{\dot H^\theta} \\
& \le \bigg\|\int_{t_1}^{t_2} (t_2-r)^{\alpha-1} (-A)^\frac\theta2 S(t_2-r) f(r) {\rm d} r\bigg\|  \\
&\quad +\bigg\|\int_0^{t_1} [(t_2-r)^{\alpha-1}-(t_1-r)^{\alpha-1}]  (-A)^\frac\theta2 S(r) f(r) {\rm d} r\bigg\|   \\
&\quad +\bigg\|\int_0^{t_1} (t_1-r)^{\alpha-1} 
[(-A)^\frac\theta2 S(t_2-r)-(-A)^\frac\theta2 S(t_1-r)] f(r) {\rm d} r\bigg\| \\
&=:I_1+I_2+I_3.
\end{align*}

Assume that $\theta\ge \rho$.
Then we have
\begin{align*}
I_1
&=\bigg\|\int_{t_1}^{t_2} (t_2-r)^{\alpha-1} 
(-A)^\frac{\theta-\rho}2 S(t_2-r) (-A)^\frac\rho 2 f(r) {\rm d} r\bigg\| \\
&\le \bigg(\int_{t_1}^{t_2} (t_2-r)^{(\alpha-1)p'} 
\|(-A)^\frac{\theta-\rho}2 S(t_2-r)\|^{p'} {\rm d} r\bigg)^\frac1{p'} \\
&\quad \times \bigg(\int_{t_1}^{t_2} \|(-A)^\frac\rho 2 f(r)\|^p {\rm d} r\bigg)^\frac1p.
\end{align*}
Since the semigroup $S(\cdot)$ is analytic, by \eqref{ana} there exists a constant $C>0$ such that 
\begin{align*}
\|(-A)^\frac{\theta-\rho}2 S(t)\|
\le C t^{-\frac{\theta-\rho}2},\quad t\in (0,T].
\end{align*}
Consequently,
\begin{align*}
I_1
&\le C\bigg(\int_0^{t_2-t_1} r^{(\alpha-1-\frac{\theta-\rho}2)p'} 
{\rm d} r\bigg)^\frac1{p'} \bigg(\int_{t_1}^{t_2} \|(-A)^\frac\rho 2 f(r)\|^p {\rm d} r\bigg)^\frac1p \\
&\le C (t_2-t_1)^{\alpha-\frac{\theta-\rho}2-\frac1p} 
\|f\|_{L^p(0,T;\dot H^\rho)}.
\end{align*}
Similarly,
\begin{align*}
I_2
&\le C \bigg(\int_0^{t_1} \frac{[(t_1-r)^{\alpha-1}-(t_2-r)^{\alpha-1}]^{p'}}  
{(t_2-r)^{\frac{(\theta-\rho)p'}2}} {\rm d} r\bigg)^\frac1{p'} 
\|f\|_{L^p(0,T;\dot H^\rho)}.
\end{align*}
Using the fact that
\begin{align*}
(b-a)^p\le b^p-a^p,\quad a\le b,\ p\ge 1,
\end{align*}
we get
\begin{align*}
I_2
&\le C \bigg(\int_0^{t_1} \frac{(t_1-r)^{(\alpha-1)p'}}
{(t_2-r)^{\frac{(\theta-\rho)p'}2}}
-(t_2-r)^{(\alpha-1-\frac{\theta-\rho}2) p'} {\rm d} r\bigg)^\frac1{p'}  
\|f\|_{L^p(0,T;\dot H^\rho)} \\
&\le C \bigg(\int_0^{t_1} (t_1-r)^{(\alpha-1-\frac{\theta-\rho}2) p'}
-(t_2-r)^{(\alpha-1-\frac{\theta-\rho}2) p'} {\rm d} r\bigg)^\frac1{p'}  
\|f\|_{L^p(0,T;\dot H^\rho)} \\
&\le C (t_2-t_1)^{\alpha-\frac{\theta-\rho}2-\frac1p} 
\|f\|_{L^p(0,T;\dot H^\rho)}.
\end{align*}
It remains to estimate $I_3$.
Note that 
\begin{align*}
(-A)^\frac{\theta-\rho}2 S(t_2-r)-(-A)^\frac{\theta-\rho}2 S(t_1-r)
=\int_{t_1-r}^{t_2-r} (-A)^{\frac{\theta-\rho}2+1} S(t) {\rm d}t.
\end{align*}
Therefore,
\begin{align*}
I_3
&\le C \int_0^{t_1} (t_1-r)^{\alpha-1} \|(-A)^\frac\rho 2 f(r)\|
\bigg(\int_{t_1-r}^{t_2-r} t^{-1-\frac{\theta-\rho}2} {\rm d}t\bigg)   
 {\rm d} r.
\end{align*}

If $\theta>\rho$, then similar arguments to estimate $I_2$ yield that
\begin{align*}
I_3
&\le C \int_0^{t_1} (t_1-r)^{\alpha-1} 
[(t_1-r)^{-\frac{\theta-\rho}2}-(t_2-r)^{-\frac{\theta-\rho}2}] 
\|(-A)^\frac\rho 2 f(r)\|  {\rm d} r \\
&\le C (t_2-t_1)^{\alpha-\frac{\theta-\rho}2-\frac1p} 
\|f\|_{L^p(0,T;\dot H^\rho)}.
\end{align*}
If $\theta=\rho$, then for any $\delta\in (0,1)$, we have
\begin{align*}
I_3
&\le C \int_0^{t_1} (t_1-r)^{\alpha-1} \|(-A)^\frac\rho 2 f(r)\|
\bigg(\int_{t_1-r}^{t_2-r} t^{-\delta} t^{-1+\delta} {\rm d}t \bigg)   
 {\rm d} r \\
&\le C \int_0^{t_1} (t_1-r)^{\alpha-1-\delta} \|(-A)^\frac\rho 2 f(r)\|
\bigg(\int_{t_1-r}^{t_2-r} t^{-1+\delta} {\rm d}t\bigg)   
 {\rm d} r \\
&\le C\delta^{-1}\int_0^{t_1} (t_1-r)^{\alpha-1-\delta} 
[(t_2-r)^\delta-(t_1-r)^\delta] \|(-A)^\frac\rho 2 f(r)\| {\rm d} r \\
&\le C \delta^{-1} (t_2-t_1)^\delta  
\bigg(\int_0^{t_1} (t_1-r)^{(\alpha-1-\delta)p'}{\rm d }r \bigg)^\frac1{p'} 
\|f\|_{L^p(0,T;\dot H^\rho)},
\end{align*}
where we use the fact that
\begin{align*}
(b-a)^p\ge b^p-a^p,\quad a\le b,\ p\le 1.
\end{align*}
Taking $\delta\in (0,\alpha-1/p)$, we obtain
\begin{align*}
I_3
&\le C (t_2-t_1)^\delta  \|f\|_{L^p(0,T;\dot H^\rho)}. 
\end{align*}

Now we assume that $\theta<\rho$.
Then we have the following estimations for the first two terms when 
$\alpha>\frac1p$:
\begin{align*}
I_1&\le C (t_2-t_1)^{\alpha-\frac1p} \|f\|_{L^p(0,T;\dot H^\rho)}, \\
I_2&\le C (t_2-t_1)^{\alpha-\frac1p} \|f\|_{L^p(0,T;\dot H^\rho)}.
\end{align*}

If $\theta\ge \rho-2\alpha+2/p$, then by \eqref{ana} there exists a constant $C>0$ such that 
\begin{align*}
\|(-A)^{1+\frac{\theta-\rho}2} S(t)\|
\le C t^{-1+\frac{\rho-\theta}2},\quad t\in [0,T].
\end{align*}
Then
\begin{align*}
I_3
&\le C \int_0^{t_1} (t_1-r)^{\alpha-1} \|(-A)^\frac\rho2 f(r)\|
\bigg(\int_{t_1-r}^{t_2-r} t^{-1+\frac{\rho-\theta}2} {\rm d}t\bigg)   
 {\rm d} r \\
 &\le C \int_0^{t_1} (t_1-r)^{\alpha-1} \|(-A)^\frac\rho2 f(r)\|
\bigg(\int_{t_1-r}^{t_2-r} t^{\frac{\rho-\theta}2-\alpha+\frac1p}
t^{-1+\alpha-\frac1p}  {\rm d}t\bigg)   
 {\rm d} r \\
&\le C \int_0^{t_1} (t_1-r)^{\frac{\rho-\theta}2-\frac1{p'}} 
[(t_2-r)^{\alpha-\frac1p}-(t_1-r)^{\alpha-\frac1p}] \|(-A)^\frac\rho 2 f(r)\| {\rm d} r \\
&\le C (t_2-t_1)^{\alpha-\frac1p} 
\bigg(\int_0^{t_1} (t_1-r)^{\frac{\rho-\theta}2p'-1} {\rm d }r \bigg)^\frac1{p'} 
\|f\|_{L^p(0,T;\dot H^\rho)} \\
&\le C (t_2-t_1)^{\alpha-\frac1p} \|f\|_{L^p(0,T;\dot H^\rho)}.
\end{align*}
Consequently, $G_\alpha$ is a bounded linear operator from $L^p(0,T; \dot H^\rho)$ to $\CC^{\alpha-1/p}([0,T]; \dot H^\theta)$ for 
$\theta\ge \rho-2\alpha+2/p$.

If $\theta<\rho-2\alpha+2/p$.
Then applying the property \eqref{ana}, we obtain the existence of a constant $C>0$ such that for any $t\in [0,T]$,
\begin{align*}
\|(-A)^{1+\frac{\theta-\rho}2} S(t)\|
\le \|(-A)^{\alpha-\frac1p-\frac{\rho-\theta}2}\|\cdot 
\|(-A)^{1-\alpha+\frac1p} S(t)\|
\le C t^{-1+\alpha-\frac1p}.
\end{align*}
Then
\begin{align*}
I_3
&\le C \int_0^{t_1} (t_1-r)^{\alpha-1} \|(-A)^\frac\rho2 f(r)\|
\bigg(\int_{t_1-r}^{t_2-r} t^{-1+\alpha-\frac1p} {\rm d}t\bigg)   
 {\rm d} r \\
&\le C \int_0^{t_1} (t_1-r)^{\alpha-1} 
[(t_2-r)^{\alpha-\frac1p}-(t_1-r)^{\alpha-\frac1p}] 
\|(-A)^\frac\rho2 f(r)\| {\rm d} r \\
&\le C (t_2-t_1)^{\alpha-\frac1p}
\bigg(\int_0^{t_1} (t_1-r)^{(\alpha-1)p'}\bigg)^\frac1{p'} 
\|f\|_{L^p(0,T;\dot H^\rho)} \\
& \le C (t_2-t_1)^{\alpha-\frac1p} \|f\|_{L^p(0,T;\dot H^\rho)}. 
\end{align*}
Thus $G_\alpha$ is a bounded linear operator from $L^p(0,T; \dot H^\rho)$ to $\CC^{\alpha-1/p} ([0,T]; \\
\dot H^\theta)$ for 
$\theta<\rho-2\alpha+2/p$.
Combining the result for $\rho-2\alpha+2/p\le \theta<\rho$, we conclude that $G_\alpha$ is a bounded linear operator from $L^p(0,T; \dot H^\rho)$ to $\CC^{\alpha-1/p} ([0,T]; \dot H^\theta)$ for 
$\theta<\rho$.
\quad \\

\begin{cor} \label{cor-OU}
Assume that $S(\cdot)$ is an analytic $\CC_0$-semigroup and there exist constants
$\alpha\in (0,1/2)$ and $\gamma\ge 0$ such that
\begin{align*}
\int_0^T r^{-2\alpha} \|S(r)\|_{\LL_2^\gamma}^2 {\rm d}r<\infty.
\end{align*}
Then for any $p\ge 1$,
\begin{align*}
W_A \in L^p(\Omega; \CC^{\delta_1} ([0,T]; \dot H^{\theta_1})) \cup
L^p(\Omega; \CC^{\delta_2}([0,T]; \dot H^{\theta_2}))
\end{align*} 
for any $\delta_1<\alpha$ with $\theta_1\in [0,\gamma]$ and $\delta_2<\alpha-(\theta_2-\gamma)/2$ with $\theta_2\in (\gamma, \gamma+2\alpha)$.
The limit case $\alpha=1/2$ is included when $(-A)^\frac\gamma2\in \LL_2^0$.
\end{cor}

{\it Proof.}
Applying Burkholder--Davis--Gundy inequality, we get
\begin{align*}
\bigg\|\int_0^t (t-r)^{-\alpha} S(t-r) {\rm d}W(r)\bigg\|_{L^p(\Omega;\dot H^\gamma)} 
\le C\bigg(\int_0^t r^{-2\alpha} \|S(r) \|_{\LL_2^\gamma}^2 {\rm d} r \bigg)^\frac12
<\infty.
\end{align*}
Then by Fubini theorem, we get
\begin{align*}
&\ee \bigg[\bigg\|\int_0^t (t-r)^{-\alpha} S(t-r) {\rm d}W(r)\bigg\|_{L^p(0,T;\dot H^\gamma)}^p  \bigg] \\
&\le \int_0^T \bigg(\int_0^t r^{-2\alpha} \|S(r)\|_{\LL_2^\sigma} {\rm d}r  \bigg)^\frac p2 {\rm d}t <\infty.
\end{align*}
This shows that $\int_0^\cdot (\cdot-r)^{-\alpha} S(\cdot-r) {\rm d}W(r)\in L^p(0,T; \dot H^\gamma)$ a.s. for any $p\ge 1$.
Now we can apply Proposition \ref{prop-hol} and obtain that 
\begin{align*}
W_A \in \CC^\delta([0,T]; \dot H^\beta) \cap 
\CC^{\alpha-\frac1p-\frac{\theta_1-\beta}2}([0,T]; \dot H^{\theta_1})\cap 
\CC^{\alpha-\frac1p}([0,T]; \dot H^{\theta_2}) 
\end{align*} 
for any $\delta<\alpha-1/p$, $\theta_1\in (\beta,\beta+2\alpha-2/p)$ and $\theta_2\in (0,\beta)$.
Applying Proposition \ref{prop-hol} with $\beta=0$, we have 
\begin{align*}
W_A \in \CC^\delta([0,T]; H) \cap 
\CC^{\alpha-\frac1p-\frac{\theta_3}2}([0,T]; \dot H^{\theta_3})
\end{align*} 
for any $\delta<\alpha-1/p$, $\theta_3\in (0,2\alpha-2/p)$.
Taking $p$ large enough, we complete the proof.

Now we assume that $(-A)^\frac\beta2\in \LL_2^0$ and $\alpha=1/2$.
Let $0\le t_1<t_2\le T$.
Then 
\begin{align*}
\ee\Big[ \|W_A(t_2)-W_A(t_1) \|_{\dot H^\beta}^2\Big]=:II_1+II_2,
\end{align*}
where 
\begin{align*}
II_1 &=\int_{t_1}^{t_2} \|(-A)^\frac\beta2 S(t_2-r)\|_{\LL_2^0}^2 
{\rm d} r, \\
II_2 &=\int_0^{t_1} \|(-A)^\frac\beta2 (S(t_2-r)-S(t_1-r)) \|_{\LL_2^0}^2 {\rm d} r.
\end{align*}
By the uniform boundedness of $\|S(\cdot)\|$, we have
\begin{align*}
II_1=C \int_{t_1}^{t_2} \|(-A)^\frac\beta2\|_{\LL_2^0}^2 {\rm d} r
\le C(t_2-t_1). 
\end{align*}
Note that 
\begin{align*}
(-A)^\frac\beta2 S(t_2-r)-(-A)^\frac\beta2 S(t_1-r)
=\int_{t_1-r}^{t_2-r} (-A)^{\frac\beta2+1} S(\rho) {\rm d}t.
\end{align*}
Then
\begin{align*}
II_2
&=\int_0^{t_1} \bigg\|\int_{t_1-r}^{t_2-r} 
(-A)^{\frac\beta2+1} S(\rho) {\rm d}t \bigg\|_{\LL_2^0}^2 {\rm d} r \\
&\le \sum_{n=1}^\infty \int_0^{t_1} \bigg(\int_{t_1-r}^{t_2-r} 
\|A S(\rho) (-A)^\frac\beta2 {\bf Q}^\frac12 e_n\| {\rm d}t \bigg)^2 {\rm d} r \\
&\le C\|(-A)^\frac\beta2\|_{\LL_2^0}^2 \int_0^{t_1} \bigg(\int_{t_1-r}^{t_2-r} 
t^{-1} {\rm d}t \bigg)^2 {\rm d} r.
\end{align*}
For any $\epsilon\in (0,1/2)$, 
\begin{align*}
II_2
&\le C\|(-A)^\frac\beta2\|_{\LL_2^0}^2 \int_0^{t_1} (t_1-r)^{-2\epsilon} \bigg(\int_{t_1-r}^{t_2-r} 
t^{\epsilon-1} {\rm d}t \bigg)^2 {\rm d} r \\
&\le C\epsilon^{-2} \|(-A)^\frac\beta2\|_{\LL_2^0}^2 \int_0^{t_1} (t_1-r)^{-2\epsilon} \big((t_2-r)^\epsilon-(t_1-r)^\epsilon\big)^2 {\rm d} r \\
&\le C\epsilon^{-2} \|(-A)^\frac\beta2\|_{\LL_2^0}^2 (t_2-t_1)^{2\epsilon}
\int_0^{t_1} (t_1-r)^{-2\epsilon} {\rm d} r \\
&\le C (t_2-t_1)^{2\epsilon}.
\end{align*}
Thus we obtain
\begin{align*}
\ee\Big[ \|W_A(t_2)-W_A(t_1) \|_{\dot H^\beta}^2\Big]
\le C (t_2-t_1)^{2\epsilon}.
\end{align*}
Since $W_A(t_2)-W_A(t_1)$ is Gaussian, we conclude that 
\begin{align*}
W_A \in  \CC^\delta ([0,T]; \dot H^\beta),\quad \delta<1/2.
\end{align*} 
Other cases that $W_A \in \CC^{\delta_1} ([0,T]; \dot H^{\theta_1})$
for any $\delta_1<(1+\beta-\theta_1)/2$ with 
$\theta_1\in (\beta, \beta+1)$ and $W_A \in \CC^{\delta_2} ([0,T]; \dot H^{\theta_2})$
for any $\delta_2<1/2$ with $\theta_2<\beta$ are analogous and we omit the details.
\quad \\

\subsection{Proof of Theorem \ref{main4}}
\label{sec4.3}

{\it Proof.}[Proof of Theorem \ref{main4}]
For the initial datum, by \eqref{ana} we get
\begin{align*}
\|S(t_2)X_0-S(t_1)X_0\|_{\theta} 
&=\|(-A)^\frac{\theta-\beta}2 (S(t_2-t_1)-{\rm Id}_H) 
(-A)^\frac\beta2 S(t_1) X_0\| \nonumber \\
&\le C|t_2-t_1|^{\frac{\beta-\theta}2\wedge 1} \|X_0\|_\beta
\end{align*}
for any $\theta\in [0,\beta)$, which combining with the proof of Theorem \ref{main2} shows that 
$S(\cdot) X_0\in \CC^{\frac{\beta-\theta}2\wedge 1}([0,T]; \dot H^\theta)$ for any $\theta\in [0,\beta]$. 

By Burkholder--Davis--Gundy inequality, we get
\begin{align*}
\|G_\alpha(t) \|_{L^p(\Omega;\dot H^{\gamma})}  
&\le \bigg(\int_0^t (t-r)^{-2\alpha} \|S(t-r) G(X(r)) \|_{L^p(\Omega;\LL_2^\gamma)}^2 {\rm d} r \bigg)^\frac12 \\
&\le \bigg(\int_0^t r^{-2\alpha} K_{G,\gamma}^2(r) {\rm d}r  \bigg)^\frac12
\bigg(1+\|X\|_{\HHH^p_\gamma}\bigg).
\end{align*}
Then by Fubini theorem, we get
\begin{align*}
& \ee \bigg[\|G_\alpha(t)\|_{L^p(0,T;\dot H^{\gamma})}^p  \bigg]  
=\int_0^T \ee\bigg[\|G_\alpha(t) \|_\gamma^p \bigg]  {\rm d}t  \\
&\le \bigg[\int_0^T \bigg(\int_0^t r^{-2\alpha} K_{G,\gamma}^2(r) {\rm d}r  \bigg)^\frac p2 {\rm d}t \bigg] \bigg(1+\|X\|_{\HHH^p_\gamma}\bigg)^p
<\infty.
\end{align*}
This shows that $G_\alpha\in L^p(\Omega;L^p(0,T; \dot H^{\gamma}))$. 

Now we can apply Proposition \ref{prop-hol} with $\rho=\gamma$.
When $\gamma=0$, we have 
\begin{align*}
S\diamond G(X) \in L^p(\Omega;\CC^\delta([0,T]; H)) \cup L^p(\Omega;\CC^{\alpha-\frac1p-\frac\theta 2}([0,T]; \dot H^\theta))
\end{align*} 
for any $\delta\in [0,\alpha-1/p)$ and $\theta\in (0,2\alpha-2/p)$.
When $\gamma>0$, we obtain 
\begin{align*}
S\diamond G(X) \in\ 
&L^p(\Omega;\CC^\delta ([0,T]; \dot H^{\gamma})) 
\cup L^p(\Omega;\CC^{\alpha-\frac1p} ([0,T]; \dot H^\theta)) \nonumber \\
& \cup L^p(\Omega;\CC^{\alpha-\frac1p+\frac{\gamma-\theta_1}2} ([0,T]; \dot H^{\theta_1}))
\end{align*} 
for any $\delta\in [0,\alpha-1/p)$, $\theta\in [0,\gamma)$ and $\theta_1\in (\gamma, \gamma+2\alpha-2/p)$.
Similar argument yields the same regularity for $S*F(X)$.
Thus we conclude the results \eqref{con-0} and \eqref{con-gamma} by combining the H\"older continuity of $S(\cdot)X_0$, $S*F(X)$ and $S\diamond G(X)$.
\quad \\

\section*{Acknowledgements}

The authors gratefully thank the anonymous referees for valuable comments and suggestions in improving the present paper.

\bibliography{bib}

\end{document}